\documentclass[11pt]{article}
\usepackage[utf8]{inputenc}
\usepackage[T1]{fontenc}
\usepackage{mathrsfs}

\usepackage{latexsym}
\usepackage{bbm}
\usepackage{amsmath}
\usepackage{amsfonts}
\usepackage{amssymb}
\usepackage{multirow}
\usepackage{latexsym}
\usepackage{overpic}
\usepackage{graphicx}
\newtheorem{theorem}{\bf Theorem}[section]
\newtheorem{lemma}[theorem]{\bf Lemma}
\newtheorem{prop}[theorem]{\bf Proposition}

\newtheorem{remark}{{\bf Remark}}[section]
\newenvironment{proof}{\noindent{\em Proof.}}{\quad \hfill$\Box$\vspace{2ex}}
\newenvironment{definition}{\noindent{\bf Definition}}{\vspace{2ex}}

\def \and {\, \mbox{\rm and}\, }

\makeatletter

\newcommand{\Rmnum}[1]{\expandafter\@slowromancap\romannumeral #1@}
\makeatother
\begin{document}

\title{On Variational Inequalities for Averages over Hypersurfaces}
\author{Xudong Nie\thanks{Corresponding author: Department of Mathematics
and Physics, Shijiazhuang Tiedao University, Shijiazhuang 050043, P.R. China. E-mail address:
{\it nxd2016@stdu.edu.cn}.}
\,\,\,\,\,\,\,\,and\,\,\,\,\,\,
Panwang Wang\thanks{Department of Mathematics
and Physics, Shijiazhuang Tiedao University, Shijiazhuang 050043, P.R. China. E-mail address:
{\it wangpw@stdu.edu.cn}.}
}
\date{}
\maketitle

%\dedicatory{}

\begin{abstract}
This paper establishes sharp variational estimates for families of averaging operators over smooth hypersurfaces in $\mathbb{R}^{n+1}$.
For a hypersurface with $k$ non-zero principal curvatures at a given point,
we prove that the $r$-variation operator $V_r(A_t f : t>0)$ is bounded on $L^p$ for all exponents satisfying $(k+1)/k < p < r(k+1)$,
with operator norm depending explicitly on a translation parameter $D$.
In the special case $k=n$ (i.e., non-vanishing Gaussian curvature),
we show that the range $p<r(n+1)$ is optimal,
thereby recovering and extending classical results for the sphere.
Furthermore, for convex hypersurfaces of finite type (in the sense of Iosevich-Sawyer),
we obtain variational bounds in the range $p_0 < p < 2r$, where $p_0$ is the contact  index of the associated hypersurface.
The proofs combine long-short variational decompositions, stationary phase analysis, vector-valued singular integrals,
and local smoothing estimates for Fourier integral operators.
These results unify and generalize earlier work on maximal functions and variational inequalities for averages over hypersurfaces.
\end{abstract}


AMS Subject Classification (2010) 42B20, 42B35

\section{Introduction}
\subsection{Motivation and background}

Let $S \subset \mathbb{R}^{n+1}$ be a smooth hypersurface and
let  $d\sigma (\overline{y})$ be a smooth measure supported in a compact subset $S_0$ of $S$.
For a Schwartz function $f$ we consider the family of averaging operators
\[
A_t f(\bar{x}) = \int_S f(\bar{x} - t\bar{y})\,d\sigma(\bar{y}),\qquad t>0,
\]
where $\overline{x}=(x,x_{n+1})\in\mathbb{R}^{n+1}$.
The study of such operators has a long history in harmonic analysis,
beginning with the seminal works of
 Stein \cite{1976-Stein-maximal-sphere} ($n>1$),  and Bourgain \cite{1986-Bourgain-maximal-plan} ($n=1$)
on spherical maximal functions.
A fundamental question,
 raised in \cite{1985-Sogge-Stein}, is to determine the range of exponents $p$ for which the maximal operator $Mf = \sup_{t>0}|A_tf|$ is bounded on $L^p$.

For general hypersurfaces the curvature plays a decisive role.
Greenleaf \cite{1981-Greenleaf-variation-k-pricipal} proved that
 if at least $k$ principal curvatures are non-zero everywhere, then the maximal operator is bounded for $p>(k+1)/k$.
 The case $k=1$ was later treated by Sogge \cite{1995-Sogge}.
  For convex hypersurfaces of finite type, Iosevich and Sawyer \cite{1997-Iosevich,1999-Iosevich-Saw} characterized $L^p$-boundedness in terms of the distance to tangent planes,
  and refined results relating the boundedness index to the oscillation index of the surface were obtained by Ikromov,
  M\"{u}ller and others \cite{2010-Ikromov-muller,2019-Buschen-part1,2019-Buschen-part2,2022-Ikromov}.

In recent years the attention has shifted to stronger variational norm estimates,
which capture oscillation along the parameter and directly imply pointwise convergence without requiring a dense class.
For a family $\{a_t(\bar{x})\}_{t>0}$ and $r\in[1,\infty)$ the $r$-variation seminorm is defined by
\[
V_r(a_t:t>0)(\bar{x}) = \sup_{N\in\mathbb{N}}\sup_{t_1<\cdots<t_N}
\Bigl(\sum_{i=1}^{N-1}|a_{t_{i+1}}(\bar{x})-a_{t_i}(\bar{x})|^r\Bigr)^{\!1/r},
\]
with the obvious modification for $r=\infty$.
The concept of $r$-variation was first used in harmonic analysis by Bourgain \cite{1989-Bourgain-variation} in the context of ergodic averages;
its roots lie in the probability and martingale theory works of L\"{e}pingle \cite{1976-lepingle-first} and Pisier-Xu \cite{1988-Pisger-xu}.
Subsequent work by Jones, Seeger and Wright \cite{2008-strongvariation} established sharp variational bounds for spherical means,
showing that $V_r(A_t:t>0)$ is bounded on $L^p$ for $(n+1)/n<p<r(n+1)$,
with the upper bound $p<r(n+1)$ being necessary. Endpoint estimates were later obtained by Beltran, Oberlin, Roncal and Seeger \cite{2022-Beltran-variation-sphere}.
In a recent paper, Nie \cite{2006-Nie-variation-plan} studied variational operators associated with planar curves of finite type,
obtaining sharp $L^p$ bounds for isotropic and non-isotropic dilations;
that work extended earlier results of Iosevich \cite{1994-Iosevich} and Li \cite{LIwenjuan} to the variational setting.
The present article continues this line of research by considering averages over higher-dimensional hypersurfaces,
allowing the curvature to vanish in a controlled way.

In addition to the aforementioned results on hypersurfaces, sharp variation-norm estimates for oscillatory integrals related to Carleson's theorem were obtained by Guo, Roos, and Yung \cite{2020-GUO-variation}, and variational inequalities for generalized spherical means were established by Li, Yang, and Zhang \cite{2024-liwenjuan-variantion}.
Dimension-free variational inequalities for averaging operators have been studied by Bourgain, Mirek, Stein and collaborators \cite{2018-Bourgain-variation},
while discrete analogues have been systematically developed by Mirek, Stein, Trojan and Zorin-Kranich
\cite{Mirek-Stein-variation-radontype,Mirek-Stein-variation-jump,Mirek-Stein-variation-bootstrap,Mirek-Stein-variation-interpolation}.
These developments highlight the growing interest in variational estimates for a wide class of averaging operators and their applications to pointwise convergence problems.

\subsection{Main results}

Let $\Phi\in C^\infty(\mathbb{R}^n)$ satisfy $\Phi(0)=0$ and $\nabla\Phi(0)=0$.
 For a fixed $D\in\mathbb{R}^+$ define the averaging operator
\[
A_tf(x,x_{n+1}) = \int_{\mathbb{R}^n} f(x-ty,\;x_{n+1}-t(D+\Phi(y)))\,\rho(y)\,dy,
\]
where $\rho$ is a smooth cutoff supported near the origin.
 Our first result concerns hypersurfaces with a prescribed number of non-zero principal curvatures.

\begin{theorem}\label{theorem-k-principal}
Suppose the hypersurface $S=\{(x,x_{n+1}):x_{n+1}=\Phi(x)+D\}$ has exactly $k$ non-zero principal curvatures at the point $(0,D)$ (with $1\le k\le n$).
Then for every $r>2$ the $r$-variation operator
satisfies
\begin{equation}\label{ineq-main-case-k}
\left\|V_{r}\left(A_tf:t>0\right)\right\|_{L^{p}(\mathbb{R}^{n+1})}
\leq C(p,r,D)\|f\|_{L^{p}(\mathbb{R}^{n+1})}
\end{equation}
 for all
\[
(k+1)/k < p < r(k+1).
\]
where the operator norm satisfies the estimate
\begin{equation}\label{C(p,r,D)-estimate}
C(p,r,D)\leq C_{\epsilon}\left(1+D^{1/(r\wedge p)+\epsilon}\right),
\end{equation}
and the growth order $1/(r\wedge p)$ in $D$ is optimal up to an arbitrarily small $\epsilon$.

Conversely, when $k=n$, if there exists a constant $C$ such that
$$
\left\|V_{r}\left(A_tf:t>0\right)\right\|_{L^{p}(\mathbb{R}^{n+1})}
\leq C\|f\|_{L^{p}(\mathbb{R}^{n+1})},
$$
then $p\leq r(n+1)$.
\end{theorem}

When $k=n$ (i.e., non-vanishing Gaussian curvature) this recovers and extends the classical spherical variational estimates of Jones-Seeger-Wright \cite{2008-strongvariation}. Moreover, the case $k<n$ was not previously known for general hypersurfaces.
Theorem \ref{theorem-k-principal} thus provides a simultaneous generalization of Greenleaf's maximal inequality
\cite{1981-Greenleaf-variation-k-pricipal} and the variational results for spheres.

Our second result deals with convex hypersurfaces of finite type, a class intensively studied by Iosevich and Sawyer \cite{1997-Iosevich}.
In the proof of this result, the sharp estimate for the dependence on the translation parameter $D$ obtained in Theorem \ref{theorem-k-principal} will play a crucial role.

\begin{theorem}\label{theorem-convex}
Let $\Phi$ be a smooth convex function and let $S=\{(x,x_{n+1}):x_{n+1}=\Phi(x)+D\}$ be its graph.
Assume $S$ has no tangent line of infinite order
at the point $(0,D)$.
Denote by $p_0$  the contact index of $S$ at (0,D), defined by
\begin{equation}\label{defi-p0}
p_0=\inf_{p}\left\{p:\int \Phi(y)^{-1/p}\rho(y)dy<\infty\right\}
\end{equation}
and suppose $p_0>2$. Then for every $r>p_0$ the estimate
$$
\left\|V_r(A_tf:t>0)\right\|_{L^p(\mathbb{R}^{n+1})}\leq C\|f\|_{L^p(\mathbb{R}^{n+1})}
$$
holds for $p_0< p<2r$.
\end{theorem}

\subsection{Strategies and novelties}
The proof of Theorem 1.1 follows a well-established strategy for variational estimates
which used in \cite{2008-strongvariation,2006-Nie-variation-plan,2020-GUO-variation,2024-liwenjuan-variantion},
which involves a long-short decomposition of the variation operator.
The novelty of our approach lies in a meticulous and sharp analysis of the dependence on the translation parameter
$D$.

\textit{Long-Short Decomposition:}
We first decompose the $r$-variation into two parts:
the ``long variation",
which captures changes across dyadic scales,
 and the ``short variation",
 which captures changes within each dyadic scale.
 The long variation is handled by straightforward decay estimates for the Fourier transform of the surface measure,
 leading to a simple bound.

\textit{Reduction to $L^2_t$ Estimates:}
 The core of the proof is the analysis of the short variation.
  We employ a powerful pointwise estimate that bounds the short variation at a dyadic scale by the product of an
$L^2_t$
  norm of the family and an
$L^2_t$
  norm of its time derivative. This cleverly reduces the problem to proving
$L^2_t$
  estimates for the averaging operators
$\mathcal{A}_{k,t}f$ and their derivatives.

\textit{Vector-Valued Fourier Multiplier and Local Smoothing:}
To obtain these $L^2_t$  estimates,
we combine several harmonic analysis tools.
 Plancherel's theorem translates the problem into frequency space.
 A detailed stationary phase analysis provides precise asymptotic for the Fourier multiplier
$m$ associated with the hypersurface.
 This, in turn, allows us to apply a vector-valued multiplier theorem and,
  crucially, local smoothing estimates for Fourier integral operators.
  These estimates are essential for handling the curvature and achieving the full range of exponents $p$.

\textit{Key Novelty - Sharp Dependence on $D$:}
A critical innovation in our argument is the careful tracking of the translation parameter
$D$  through every step-from the stationary phase expansion, through the multiplier estimates, to the final
$L^p$ bounds.
This level of detail is not present in previous works.
The resulting estimate
$C(p,r,D)\lesssim 1+D^{\frac{1}{r\wedge p}+\epsilon}$
is not only optimal in its growth order,
as we prove via counterexamples,
but also becomes an indispensable tool for extending our results to more complex surfaces.
This sharp $D$-dependence is what allows us to apply a scaling argument to convex finite-type hypersurfaces in Theorem \ref{theorem-convex},
where the translation parameter is effectively scaled up to $2^kD$.

\subsection{Organisation of the paper}
The paper is structured as follows.

\textbf{Section 2: Preliminaries.}
We introduce the key notational conventions used throughout the paper.
We then define the long and short variational decompositions and recall a fundamental estimate for bounding the short variation in terms of
$L^2_t$ norms.

\textbf{Section 3: The case $k=n$.} This section is dedicated to the core proof of Theorem \ref{theorem-k-principal}
 for hypersurfaces with non-vanishing Gaussian curvature.
 It is divided into several subsections that handle the long variation, the short variation, and the necessary
$L^2$, $L^1$, and $L^p$ estimates for the decomposed operators.
It concludes with the proof of the necessary condition $p\leq r(n+1)$.

\textbf{Section 4: The case $k<n$.}
Here, we extend the result of Section 3 to hypersurfaces with only $k$ non-zero principal curvatures.
The main tool is a rotation argument (Lemma 4.1),
which averages the operator over rotations in the ``flat" directions,
reducing the problem to the previously established full-curvature case.

\textbf{Section 5: Convex finite-type surfaces. }
We prove Theorem \ref{theorem-convex} for convex hypersurfaces of finite type.
This is achieved by using a non-isotropic dyadic decomposition and a powerful scaling argument.
 The scaling property, combined with the sharp $D$-dependence from Theorem \ref{theorem-k-principal},
  yields the desired variational bounds.

\textbf{Section 6: Appendix. }
The appendix contains the proof of the auxiliary estimates for the Poisson-like operator $S_t$,
which were used in Section 3 to simplify the main operator.

\section{Preliminaries}
This section collects the essential definitions, notations, and preliminary lemmas that will be used throughout the paper.

\subsection{Notation}
For any two non-negative quantities $x$, $y$ we  write $x\lesssim y$ (or $x\gtrsim y$) to denote that
there exists an absolute constant $C>0$ such that
 $x\leq Cy$ (or $x\geq Cy$).
If the constant depends on a parameter,  e.g., $\epsilon$,
we write $x\lesssim_{\epsilon}y$.

For a sequence $\{a_k\}_{k\in\mathbb{Z}}$, we denote its $l^2$-norm by $\|a_k\|_{l^r_k}:=\left(\sum_{k\in\mathbb{Z}}|a_k|^r\right)^{1/r}$.
For a family of functions indexed by a continuous parameter $t>0$,
we use the notation
$\|a_t\|_{L^r_t}:=\left(\int_0^{\infty}|a_t|^r\frac{dt}{t}\right)^{1/r}$.

We fix a smooth cutoff function $\psi:\mathbb{R}\rightarrow[0,1]$  supported in $[1/2, 2]$ such that
$\sum_{j\in\mathbb{Z}}\psi_j(t)=1$ for all $t>0$.
We set $\psi_{j}(t)=\psi\left(2^{-j}t\right)$.
This partition of unity allows us to define dyadic frequency projections on
$\mathbb{R}^{n+1}$.
For each $j\in\mathbb{Z}$,
let $P_j$ be the Fourier multiplier operator defined by
$$\widehat{P_jf}(\xi,\xi_{n+1})=\psi_{j}(\xi,\xi_{n+1})\widehat{f}(\xi,\xi_{n+1}).$$

\subsection{The Long-Short decomposition}
A key tool in proving variational estimates is the decomposition of the variation into ``long" and ``short" parts.
For a family of operators $\{A_t\}_{t>0}$,
define the short variation on the dyadic interval $[2^j,2^{j+1}]$ as
$$
V_r^j(Af)(\overline{x}):=V_r\left(A_tf(\overline{x}):t\in[2^j,2^{j+1}]\right).
$$
The total short variation is then given by the $l^r$-sum of these dyadic contributions:
$$
S_r(A_tf(\overline{x}):t>0):=\left(\sum_{j\in\mathbb{Z}}|V_r^j(Af)(\overline{x})|^r\right)^{1/r}.
$$

The following lemma justifies this decomposition.

\begin{lemma}
For any $r\in[1,\infty)$, there  exists a constant $C_r>0$ such that
\begin{equation}\label{ineq-long-short}
V_r\left(A_tf(\overline{x}):t>0\right)\leq C_r S_r(A_tf(\overline{x}):t>0)+C_r V_r\left(A_{2^j}f(\overline{x}):j\in\mathbb{Z}\right).
\end{equation}
where $V_r\left(A_{2^j}f(\overline{x}):j\in\mathbb{Z}\right)$
is the $r$-variation of the sequence of operators sampled at the dyadic points $t=2^j$.
\end{lemma}

The first term on the right-hand side captures the oscillation of
$A_tf$ within each dyadic scale,
while the second term captures the longer-range variation between scales.
This decomposition reduces the problem to estimating these two separate components.

\subsection{A pointwise estimate for the short variation}
A crucial estimate for controlling the short variation is provided by the following proposition.
It bounds the short variation at a given scale by the $L^2_t$
  norms of the function and its time derivative.

\begin{prop}\label{prop-variation-estimate}
For any $r>1$,
let $r'$  be its conjugate exponent (i.e., $\frac{1}{r}+\frac{1}{r'}=1$).
Then there exists a constant
$C_r>0$ such that for any differentiable familly $\{a_t\}_{t>0}$,
\begin{equation}\label{ineq-sr-estimate}
S_r(a_{t}:t>0)\leq C_r \|a_t\|_{L^r_t}^{1/r}\|t\cdot a'_t\|_{L^r_t}^{1/r'}.
\end{equation}
\end{prop}

A more convenient form for our purposes,
which follows directly from H\"{o}lder's inequality,
is the following pointwise estimate for the dyadic short variation:
\begin{equation*}
S_2(A_{t}f(\overline{x}):t>0)\lesssim \left(\int_{0}^{\infty}
|A_{t}f(\overline{x})|^2\frac{dt}{t}\right)^{1/4}
\left(\int_{0}^{\infty}
\left|t\frac{d}{dt}A_{t}f(\overline{x})\right|^2\frac{dt}{t}\right)^{1/4}.
\end{equation*}
This inequality is fundamental as it transforms the problem of estimating a complicated variational norm into the more tractable problem of estimating
$L^2_t$
  norms of the operator and its derivative. The majority of the technical work in Section 3 will be dedicated to establishing
$L^p$
  bounds for these
$L^2$
  norms.

  \section{The proof of Theorem \ref{theorem-k-principal}: the case $k=n$}
In this section we consider the case $k=n$ in Theorem \ref{theorem-k-principal}.
Since the  norm $\|\cdot\|_{L^p\rightarrow L^p}$ of operator $V_{r}({A_t}:t>0)$ remains unchanged under orthogonal transformation, we can assume
$$
\Phi(x)=\sum_{i=1}^{n}\lambda_ix_i^2+R(x)
$$
where $\lambda_i\neq 0$ and $D^{\alpha}R(0)=0$ for $|\alpha|\leq2$.
Without loss of generality, assume $\int_{\mathbb{R}^n}\rho(x)dx=1$.
Thus, if we let
 $
\mu(\xi,\xi_{n+1})=\int_{\mathbb{R}^n}e^{-i\left(x\cdot\xi+\xi_{n+1}(D+\Phi(x))\right) }
\rho(x)dx,
$
then we obtain $\mu(0,0)=1$.

In this section we shall fix $r\in(2,\infty)$.
Denote
$$\mathcal{F}(S_tf)(\xi,\xi_{n+1})=e^{-t|(\xi,\xi_{n+1})|}e^{-itD\xi_{n+1}}\widehat{f}(\xi,\xi_{n+1})$$
and
$$
\mathcal{F}(\mathcal{A}_tf)(\xi,\xi_{n+1})=m(t\xi,t\xi_{n+1})\widehat{f}(\xi,\xi_{n+1})
$$
where $m(\xi,\xi_{n+1})=\mu(\xi,\xi_{n+1})-e^{-|(\xi,\xi_{n+1})|}e^{-iD\xi_{n+1}}$.
For every $f\in L^p(\mathbb{R}^{n+1})$, we obtain
\begin{equation}\label{ineq-zhuanhua}
\|V_r(A_tf:t>0)\|_{L^p}\leq C_r\|V_r(S_tf:t>0)\|_{L^p}+C_r\|V_r(\mathcal{A}_tf:t>0)\|_{L^p}.
\end{equation}
In Appendix, we shall prove the estimate
$$
\|V_r(S_tf:t>0)\|_{L^p}\leq C(p,r,D)\|f\|_{L^p}
$$
for $p>1$, where $C(p,r,D)$ satisfies (\ref{C(p,r,D)-estimate}).
Therefore, it remains to obtain $L^p$ bounds for the second term in (\ref{ineq-zhuanhua}).
For this purpose, we  use (\ref{ineq-long-short}).
Indeed, we shall show that,
\begin{equation}\label{ineq-case-long}
\|V_r(\mathcal{A}_{2^j}f:j\in\mathbb{Z})\|_{L^p}\leq C\|f\|_{L^p}
\end{equation}
for $p>1$
and
\begin{equation}\label{ineq-case-short}
\|S_r(\mathcal{A}_{t}f:t>0)\|_{L^p}\leq C(p,r,D)\|f\|_{L^p},
\end{equation}
for $(n+1)/n<p<r(n+1)$.

\subsection{Estimates for long variations: proof of (\ref{ineq-case-long})}
Recall that
$
V_r(\mathcal{A}_{2^j}f:j\in\mathbb{Z})\leq \|\mathcal{A}_{2^j}f\|_{l^r_j}\leq \|\mathcal{A}_{2^j}f\|_{l^2_j}.
$
To prove (\ref{ineq-case-long}) we shall show the following bounds
\begin{equation}\label{ineq-long-square}
\|\mathcal{A}_{2^j}f\|_{L^p(l^2_j)}
\leq C\|f\|_{L^p},
\end{equation}
for $p>1$.
In order to justify (\ref{ineq-long-square}), we   shall use the Littlewood-Paley operator $P_j$.
Indeed, observe that
\begin{align}\label{ineq-long-shoulian}
&\|\mathcal{A}_{2^j}f\|_{L^p(l^2_j)}
=\|\sum_{k\in\mathbb{Z}}\mathcal{A}_{2^j}P_{-(j+k)}f\|_{L^p(l^2_j)}
\leq \sum_{k\in\mathbb{Z}}\|\mathcal{A}_{2^j}P_{-(j+k)}f\|_{L^p(l^2_j)}\nonumber\\
\leq& C\sum_{k\in\mathbb{Z}}2^{-\delta_{p}|k|}\|f\|_{L^p}\leq C\|f\|_{L^p}.
\end{align}
The proof of the last but one inequality in (\ref{ineq-long-shoulian}) follows from
\begin{equation}\label{ineq-long-22}
\|\mathcal{A}_{2^j}P_{-(j+k)}f\|_{L^2(l^2_j)}\leq C2^{-\delta|k|}\|f\|_{L^2},
\end{equation}
and
\begin{equation}\label{ineq-long-pp}
\|\mathcal{A}_{2^j}P_{-(j+k)}f\|_{L^p(l^2_j)}\leq C\|f\|_{L^p}.
\end{equation}

To prove (\ref{ineq-long-22}),
we firstly show that
\begin{equation}\label{ineq-m-pr}
|m(\xi,\xi_{n+1})|\leq C|(\xi,\xi_{n+1})|,\,\,\,\,\,\,\,\,\,\, |m(\xi,\xi_{n+1})|\leq C\left(1+|(\xi,\xi_{n+1})|\right)^{-\frac{n}{2}}
\,\,\,\,\,\,\,\,\,\,
\end{equation}
The first estimate in (\ref{ineq-m-pr})  follows from
\begin{equation*}
\left|\mu(\xi,\xi_{n+1})-e^{-iD\xi_{n+1}}\right|=\left|\int\left(e^{- i(x\cdot \xi+\xi_{n+1}\Phi(x))}-1\right)e^{-iD\xi_{n+1}}\rho(x)dx\right|\leq C|(\xi,\xi_{n+1})|
\end{equation*}
and
\begin{equation*}
\left|e^{-|(\xi,\xi_{n+1})|}e^{-iD\xi_{n+1}}-e^{-iD\xi_{n+1}}\right|\leq C|(\xi,\xi_{n+1})|.
\end{equation*}
The second estimate in (\ref{ineq-m-pr})  follows from the property of $e^{-|(\xi,\xi_{n+1})|}$ and
$$
|\mu(\xi,\xi_{n+1})|\leq C\left(1+|(\xi,\xi_{n+1})|\right)^{-\frac{n}{2}},
$$
which can be obtained by the method of stationary phase.
Therefore, using  Plancherel's theorem we have
\begin{align}\label{001}
\|\mathcal{A}_{2^j}P_{-(j+k)}f\|_{L^2(l^2_j)}= &\|m(2^j\cdot)\psi(2^{j+k}\cdot)\widehat{f}(\cdot)\|_{L^2(l^2_j)}\nonumber\\
\leq&\sup_{(\cdot)\in\mathbb{R}^{n+1}}\|m(2^j\cdot)\psi(2^{j+k}\cdot)\|_{l^2_j}\|f\|_{L^2}.
\end{align}
Fix $(\cdot)\in\mathbb{R}^{n+1}$ and $k\in\mathbb{Z}$, (\ref{ineq-m-pr}) implies
\begin{align}\label{002}
\|m(2^j\cdot)\psi(2^{j+k}\cdot)\|_{l^2_j}&\leq C\min\{2^{-k},(1+2^{-k})^{-n/2}\}\nonumber\\
&\leq C\min\{2^{-k},(1+2^{-k})^{-1/2}\}\leq C2^{-\frac12|k|}.
\end{align}
(\ref{001}) and (\ref{002}) yield (\ref{ineq-long-22}).

To prove (\ref{ineq-long-pp}), it is enough to  show
\begin{equation}\label{ineq-long-xiangliang-pp}
\|\mathcal{A}_{2^j}g_j\|_{L^p(l^2_j)}\leq C\|g_j\|_{L^p(l^2_j)}.
\end{equation}
Thus in view of (\ref{ineq-long-xiangliang-pp}) and the Littlewood-Paley theorem, we conclude
\begin{equation*}
\|\mathcal{A}_{2^j}P_{-(j+k)}f\|_{L^p(l^2_j)}\leq C\|P_{-(j+k)}f\|_{L^p(l^2_j)}\leq C\|f\|_{L^p}
\end{equation*}
which proves (\ref{ineq-long-pp}).

Now, we estimate (\ref{ineq-long-xiangliang-pp}).
This method comes from \cite{1986-Duoandikoetxea}.
By duality, it suffices to consider the case $p\geq2$.
Choose $q=(p/2)'$ and let $\mathcal{H}$ denote the tempered distribution whose Fourier transform is $m$
and $\mathcal{H}^{(k)}(\cdot)=2^{-k(n+1)}\mathcal{H}(2^{-k}\cdot)$.
Then there exists $g\in L^q$ with unit norm such that
\begin{align}\label{ineq-dui}
\left\|\mathcal{A}_{2^k}g_k\right\|_{L^p(l_k^2)}^2&=\int\sum_{k\in\mathbb{Z}}|\mathcal{A}_{2^k}g_k|^2g=\int\sum_{k\in\mathbb{Z}}|\mathcal{H}^{(k)}*g_k|^2g\nonumber\\
&\leq \sum_{k\in\mathbb{Z}}\int|\mathcal{H}^{(k)}|*|g_k|^2g\leq \int\sum_{k\in\mathbb{Z}}|g_k|^2|\mathcal{H}^{(k)}|*g\nonumber\\
&\leq \|g_k\|_{L^{p}(l_k^2)}^2\|\sup_{k\in\mathbb{Z}}|\mathcal{H}^{(k)}*g|\|_{L^q}.
\end{align}
Recall that (\ref{ineq-long-22}) implies, for $q=2$,
\begin{equation}\label{ineq-maximal}
\left\|\sup_{k}|\mathcal{H}^{(k)}*g|\right\|_{L^q}\leq \left\|\mathcal{A}_{2^k}g\right\|_{L^q(l^2_k)}\leq C\|g\|_{L^q}.
\end{equation}
(\ref{ineq-dui}) and (\ref{ineq-maximal}) yield (\ref{ineq-long-xiangliang-pp}) holds for $4/3<p<4$.
This result implies (\ref{ineq-maximal}) holds for $4/3<p<4$ in turn.
Thus we can obtain a series of intervals
$$
\left(\frac43,4\right)\subset\left(\frac{8}{7},8\right)\subset\cdots\subset\left(\frac{2j}{2j-1},2j\right)\subset
\cdots
$$
which ensure (\ref{ineq-long-xiangliang-pp}) holds for $p\in\left(\frac{2j}{2j-1},2j\right)$.
It is easy to show $\lim_{j\rightarrow\infty}2j=\infty$.
Thus (\ref{ineq-long-xiangliang-pp}) holds for $p>1$.
The proof of (\ref{ineq-case-long}) is completed.

\subsection{Estimates for short variations: proof of (\ref{ineq-case-short})}
This subsection is devoted to proving (\ref{ineq-case-short}).
We firstly present an elementary decomposition of $\mathcal{A}_tf$.
Using the partition of unity $\sum_{k\in\mathbb{Z}}\psi({2^{-k}\cdot})=1$,
we get
\begin{equation}\label{ineq-unity}
\mathcal{A}_{t}f(\overline{x}):=\sum_{k\in\mathbb{Z}}\mathcal{A}_{k,t}f(\overline{x})
\end{equation}
where
\begin{equation*}
\mathcal{F}(\mathcal{A}_{k,t}f)(\xi,\xi_{n+1})=m(t\xi,t\xi_{n+1})\psi(2^{-k}t\xi,2^{-k}t\xi_{n+1})\widehat{f}(\xi,\xi_{n+1}).
\end{equation*}
Using interpolation, it suffice to prove  Lemmas \ref{le-0101}, \ref{le-0102}, \ref{le-0103} and \ref{le-0104} below.
\begin{lemma}\label{le-0101}
For $k\in\mathbb{Z}$, we have
\begin{equation}\label{ineq-short-22}
\left\|V_2^{j}(\mathcal{A}_kf)\right\|_{L^2(l^2_j)}\leq C(1+D^{1/2})\min\left\{ 2^{-k\cdot\frac{n-1}{2}},2^k\right\}\|f\|_{L^2}.
\end{equation}
\end{lemma}
\begin{lemma}\label{le-0102}
For $k\in\mathbb{Z}$ and an arbitrary positive number $\epsilon$, there exists a constant $C_{\epsilon}$ such that
\begin{equation}\label{ineq-short-11}
\left\|V_2^{j}(\mathcal{A}_kf)\right\|_{L^1(l^2_j)}\leq C_{\epsilon}(1+D^{1+\epsilon})\max\left\{2^{k+\epsilon},2^{-k\epsilon}\right\}\|f\|_{H^1}.
\end{equation}
\end{lemma}
\begin{lemma}\label{le-0103}
Suppose that $n\geq2$ and $p>2+\frac{4}{n}$, or $n=1$, $p\geq 4$. Then for $k\in\mathbb{Z}$, we obtain
\begin{equation}\label{ineq-short-pp}
\left\|V_2^{j}(\mathcal{A}_kf)\right\|_{L^p(l^2_j)}\leq C_{\epsilon}(1+D^{1/2}) \min\left\{2^{-k\left(\frac{n+1}{p}-\frac12-\epsilon\right)},2^{k\epsilon}\right\}\|f\|_{L^p}.
\end{equation}
\end{lemma}
\begin{lemma}\label{le-0104}
For $n=1$ and $k\in\mathbb{Z}$, we have
\begin{equation}\label{ineq-short-44}
\left\|V_4^{j}(\mathcal{A}_kf)\right\|_{L^4(l^4_j)}\leq C_{\epsilon}(1+D^{1/4}) \min\left\{2^{-k\left(\frac{1}{8}-\epsilon\right)},1\right\}\|f\|_{L^4}.
\end{equation}
\end{lemma}

Now we can prove (\ref{ineq-case-short}), assuming the validity of lemmas above.

\textit{Proof of (\ref{ineq-case-short})}:
We first consider the case $n\geq2$.
Interpolating (\ref{ineq-short-22}) and (\ref{ineq-short-11}), we obtain that there exists a positive constant $\delta_{p}$ such that
\begin{equation}\label{080901}
\left\|V_2^{j}(\mathcal{A}_kf)\right\|_{L^p(l^2_j)}\leq C_{\epsilon}(1+D^{1/p+\epsilon})2^{-\delta_p|k|}\|f\|_{L^p},
\end{equation}
for  $(n+1)/n<p\leq2$.
Interpolating (\ref{ineq-short-22}) and (\ref{ineq-short-pp}), we get there exists a constant $\delta_{p}$ such that
\begin{equation}\label{080902}
\left\|V_2^{j}(\mathcal{A}_kf)\right\|_{L^p(l^2_j)}\leq C(1+D^{1/2})2^{-\delta_p|k|}\|f\|_{L^p},
\end{equation}
for $2< p<2(n+1)$.
Thus (\ref{080901}) and (\ref{080902}) yield, for $(n+1)/n< p<2(n+1)$,
\begin{equation}\label{080903}
\left\|V_2^{j}(\mathcal{A}f)\right\|_{L^p(l^2_j)}\leq C(p,2,D)\|f\|_{L^p}.
\end{equation}
(\ref{ineq-case-short}) can be deduced since  $\|\cdot\|_{l^{p_1}}\leq\|\cdot\|_{l^{p_2}}$ for $p_2<p_1$ and
\begin{equation}\label{081701}
\left\|V_{\infty}^{j}(\mathcal{A}f)\right\|_{L^{\infty}(l^{\infty}_j)}\leq C\|f\|_{L^{\infty}}.
\end{equation}

Now we consider  the case $n=1$.
By interpolation, (\ref{ineq-short-22}) and (\ref{ineq-short-pp}) imply
\begin{equation}\label{080911}
\left\|V_2^{j}(\mathcal{A}_kf)\right\|_{L^p(l^2_j)}\leq C_{\epsilon}(1+D^{1/2})2^{\epsilon k}\|f\|_{L^p}
\end{equation}
for $2< p\leq4$.
Using interpolation again, (\ref{ineq-short-44}) and (\ref{080911}) yield
\begin{equation}\label{081702}
\left\|V_{r}^{j}(\mathcal{A}_kf)\right\|_{L^p(l^{2}_j)}\leq C_{\epsilon}(1+D^{1/r})2^{-\delta_{r,p}|k|}\|f\|_{L^p}
\end{equation}
for $2<r\leq 4$ and $r<p\leq 4$.
Thus (\ref{081702}) and (\ref{081701}) yield
\begin{equation}\label{080912}
\left\|V_r^{j}(\mathcal{A}f)\right\|_{L^p(l^r_j)}\leq C(1+D^{1/r})\|f\|_{L^p},
\end{equation}
for $2< p<2r$.
Then (\ref{ineq-case-short}) follows from (\ref{080912}) and the proof is complete.

\subsection{Proof of Lemma \ref{le-0101}.}
 Using (\ref{ineq-sr-estimate}) and H\"{o}lder inequality, we get
\begin{equation}\label{ineq-081501}
\left\|V_2^{j}(\mathcal{A}_kf)\right\|_{L^2(l^2_j)}\leq C\left\|\mathcal{A}_{k,t}f\right\|_{L^2(L^2_t)}^{1/2}\left\|t\frac{d}{dt}\mathcal{A}_{k,t}f\right\|_{L^2(L^2_t)}^{1/2}.
\end{equation}
Thus, we obtain from (\ref{ineq-081501}) via Plancherel's theorem that
\begin{align} \label{ineq-020801}
\left\|V_2^{j}(\mathcal{A}_kf)\right\|_{L^2(l^2_j)}&\leq C\left\|m_k(t\cdot)\widehat{f}(\cdot)\right\|_{L^2(L^2_t)}^{1/2}\left\|(t\cdot)\nabla m_k(t\cdot)\widehat{f}(\cdot)\right\|_{L^2(L^2_t)}^{1/2}\nonumber\\
&\leq C\sup_{|\cdot|=1}\|m_{k}(t\cdot)\|_{L^2_t}^{1/2}\sup_{|\cdot|=1}\|(t\cdot)\nabla m_k(t\cdot)\|_{L^2_t}^{1/2}\|f\|_{L^2},
\end{align}
where $m_{k}(\cdot)=m(\cdot)\psi(2^{-k}\cdot)$.
Thus it is important to estimate $\|m_k\|_{L^{\infty}}$ and $\|(\cdot)\nabla m_k(\cdot)\|_{L^{\infty}}$,
since both $m_k$ and $(\cdot)\nabla m_k(\cdot)$ are supported where $|\cdot|\sim 2^k$.

By (\ref{ineq-m-pr}), we can deduce that
\begin{equation}\label{0817a}
\|m_k\|_{L^{\infty}}\leq C\min\{2^{-(nk)/2},2^k\}.
\end{equation}
Rewrite $m(\xi,\xi_{n+1})=e^{-iD\xi_{n+1}}\left(\int e^{-i(x\cdot\xi+\xi_{n+1}\Phi(x))}\rho(x)dx-e^{-|(\xi,\xi_{n+1})|}\right)$.
Using the method of stationary phase for $|(\xi,\xi_{n+1})|$ near $\infty$ and computing directly for
$|(\xi,\xi_{n+1})|$ near zero, we obtain
\begin{equation}\label{0817b}
|(\xi,\xi_{n+1})\cdot\nabla m(\xi,\xi_{n+1})|\leq C(1+D)\min\left\{|(\xi,\xi_{n+1})|,|(\xi,\xi_{n+1})|^{-\frac{n}{2}+1}\right\}.
\end{equation}
It follows from (\ref{0817a}) and (\ref{0817b}) that
\begin{equation}\label{ineq-081502}
\|m_{k}(t\cdot)\|_{L^2_t}\leq C\min\{2^{-(nk)/2},2^k\}\,\,\,\,\,\,\mbox{and}\,\,\,\,\,\,\|(t\cdot)(\nabla m_k)(t\cdot)\|_{L^2_t}\leq C(1+D)\min\{2^{-(nk)/2+k},2^k\}.
\end{equation}
(\ref{ineq-020801}) and (\ref{ineq-081502}) yield (\ref{ineq-short-22}).
This completes the proof.

\subsection{Proof of Lemma \ref{le-0102}.}
Denote $\widehat{K_k}(\xi,\xi_{n+1})=\mu(\xi,\xi_{n+1})\psi(2^{-k}\xi,2^{-k}\xi_{n+1})$
and $K_{k,t}(\overline{x})=t^{-{(n+1)}}K_{k}(\frac{\overline{x}}{t})$.
Moreover, let
$\widehat{G_k}(\xi,\xi_{n+1})=e^{-|\overline{\xi}|-iD\xi_{n+1}}\psi(2^{-k}\xi,2^{-k}\xi_{n+1})$
and $H_{k}=K_k-G_k$.
To achieve (\ref{ineq-short-11}), we shall use vector-valued singular integral theory
for the convolution operators $\{\mathcal{A}_{k,t}f=H_{k,t}*f\}_{t>0}$.
Here we take the vector valued kernel  $H_k=\{H_{k,t}\}_{t>0}$
and denote the norm by $S_2(H_{k,t}:t>0)$.
Recalling that
(\ref{ineq-short-22}) gives
$$
\left\|S_2\left(H_{k,t}*f:t>0\right)\right\|_{L^2}
\leq C_{\epsilon}(1+D^{1+\epsilon})\max\left\{2^{k+\epsilon},2^{-k\epsilon}\right\}\|f\|_{L^2}.
$$
By Calder\'{o}n-Zygmund theory,
it remains  to prove
\begin{equation}\label{081801-0}
I_{H}(\bar{y}):=\int_{|\bar{x}|>2|\bar{y}|} S_2\left(H_{k,t}(\bar{x}-\bar{y})-H_{k,t}(\bar{x}):t>0\right)d\bar{x}\leq  C_{\epsilon}(1+D^{1+\epsilon})\max\left\{2^{k+\epsilon},2^{-k\epsilon}\right\}.
\end{equation}
Without loss of generality, we only prove
\begin{equation}\label{081801}
I(\bar{y}):=\int_{|\bar{x}|>2|\bar{y}|} S_2\left(K_{k,t}(\bar{x}-\bar{y})-K_{k,t}(\bar{x}):t>0\right)d\bar{x}\leq  C_{\epsilon}(1+D^{1+\epsilon})\max\left\{2^{k+\epsilon},2^{-k\epsilon}\right\}.
\end{equation}

Using $S_2(a_t:t>0)\leq C\|ta_t'\|_{L^1_t}$,
it can be deduced  that
\begin{align*}\label{081802}
I(\bar{y})\leq &C\int_{|\bar{x}|>2|\bar{y}|}\int_0^{\infty}t\left|\frac{d}{dt}\left(K_{k,t}(\bar{x}-\bar{y})-K_{k,t}(\bar{x})\right)\right|\frac{dt}{t}d\bar{x}\nonumber\\
\leq &C(n+1)\int_0^{\infty}\int_{|\bar{x}|>2|\bar{y}|}\left|K_{k,t}(\bar{x}-\bar{y})-K_{k,t}(\bar{x})\right|d\bar{x}\frac{dt}{t}\\
&+\int_0^{\infty}\int_{|\bar{x}|>2|\bar{y}|}\left|\widetilde{K}_{k,t}(\bar{x}-\bar{y})-\widetilde{K}_{k,t}(\bar{x})\right|d\bar{x}\frac{dt}{t}
\end{align*}
where
$$
\widetilde{K}_{k,t}(\bar{x})=\bar{x}\cdot \nabla_{\bar{x}}K_{k,t}(\bar{x}).
$$
We only need to estimate the last double integrals.
In doing so, without loss of generality,
we shall show
\begin{equation}\label{ineq-0209-01}
\int_0^{1}\int_{|\bar{x}|>2|\overline{y}|}\left|\widetilde{K}_{k,t}(\bar{x}-\bar{y})-\widetilde{K}_{k,t}(\bar{x})\right|d\bar{x}\frac{dt}{t}
\leq  C_{\epsilon}(1+D^{1+\epsilon})\max\left\{2^{k+\epsilon},2^{-k\epsilon}\right\}
\end{equation}
and
\begin{equation}\label{ineq-0209-02}
\int_1^{\infty}\int_{|\bar{x}|>2|\overline{y}|}\left|\widetilde{K}_{k,t}(\bar{x}-\bar{y})-\widetilde{K}_{k,t}(\bar{x})\right|d\bar{x}\frac{dt}{t}
\leq  C_{\epsilon}(1+D^{1+\epsilon})\max\left\{2^{k+\epsilon},2^{-k\epsilon}\right\}.
\end{equation}

To achieve (\ref{ineq-0209-01}), it is no restriction to  assume that $|\bar{y}|=1$.
Thus we obtain
\begin{align}\label{ineq-020904}
&\int_0^{1}\int_{|\bar{x}|>2}\left|\widetilde{K}_{k,t}(\bar{x}-\bar{y})-\widetilde{K}_{k,t}(\bar{x})\right|d\bar{x}\frac{dt}{t}\nonumber\\
&=\int_0^{1}\int_{|\bar{x}|>2/t}\left|\widetilde{K}_{k}\left(\bar{x}-\frac{\bar{y}}{t}\right)-\widetilde{K}_{k}(\bar{x})\right|d\bar{x}\frac{dt}{t}\nonumber\\
&\leq2\int_0^{1}\int_{|\bar{x}|>1/t}\left|\widetilde{K}_{k}(\bar{x})\right|d\bar{x}\frac{dt}{t}=2\int_{|\bar{x}|>1}\left|\widetilde{K}_{k}(\bar{x})\right|
\log|\bar{x}|d\bar{x}.
\end{align}

To continue our proof, let us analyze the expression of $\widetilde{K}_{k}$.
Denote $\phi(\overline{x})=\check{\psi}(\overline{x})$ and $\phi_i(\overline{x})=\frac{\partial}{\partial x_i}\phi(\overline{x})$.
Then we have
\begin{align}\label{ineq-0210-01}
K_{k}(\overline{x})
%&=\frac{1}{(2\pi)^{n+1}}\int e^{i\bar{x}\cdot\bar{\xi}}\mu(\bar{\xi})\psi(2^{-k}\bar{\xi})d\bar{\xi}\\
=\check{\mu}*\phi_{(2^k)}(\bar{x})=2^{k(n+1)}\int\phi(2^k(x-y,x_{n+1}-(D+\Phi(y))))\rho(y)dy,
\end{align}
where  $\phi_{(2^k)}(\overline{x})=\phi(2^k\overline{x})$.
It follows that
\begin{equation}\label{ineq-020903}
\widetilde{K}_{k}(\bar{x})
%&=\sum_{i=1}^{n+1}x_i\frac{\partial}{\partial x_i}K_{k}(\bar{x})\\
=\sum_{i=1}^{n+1}2^{k(n+2)}x_i\int\phi_i(2^k(x-y,x_{n+1}-(D+\Phi(y))))\rho(y)dy.
\end{equation}

Using (\ref{ineq-020904}) and (\ref{ineq-020903}), we have
\begin{align*}
&\int_{|\bar{x}|>1}\left|\widetilde{K}_{k}(\bar{x})\right|
\log|\bar{x}|d\bar{x}\\
%&\leq \int_{|\bar{x}|>1}\int\sum_{i=1}^{n+1}2^{k(n+2)}\left|x_i\phi_i(2^k(x-y,x_{n+1}-(D+\Phi(y))))\right|\rho(y)
%\log|\bar{x}|dy d\bar{x}\\
&= \int\int_{|\bar{x}|>1}\sum_{i=1}^{n+1}2^{k(n+2)}\left|x_i\phi_i(2^k(x-y,x_{n+1}-(D+\Phi(y))))\right|
\log|\bar{x}|d\bar{x}\rho(y)dy\\
&= \int\int_{|\bar{x}|>2^k}\sum_{i=1}^{n+1}\left|x_i\phi_i(x-2^ky,x_{n+1}-2^k(D+\Phi(y)))\right|
\log|2^{-k}\bar{x}|d\bar{x}\rho(y)dy.
\end{align*}
By inequality $\log(|x|)\leq C_{\epsilon}|x|^{\epsilon}$ for $|x|\geq 1$,
it can be deduced that
\begin{align*}
\int_{|\bar{x}|>1}\left|\widetilde{K}_{k}(\bar{x})\right|
\log|\bar{x}|d\bar{x}
&\leq C_{\epsilon}\int\int\sum_{i=1}^{n+1}\left|x_i(2^{-k}|\overline{x}|)^{\epsilon}\phi_i(x-2^ky,x_{n+1}-2^k(D+\Phi(y)))\right|
d\bar{x}\rho(y)dy\\
&\leq C_{\epsilon}(1+D^{1+\epsilon})\max\left\{2^{k},2^{-k\epsilon}\right\}.
\end{align*}

Now, we prove (\ref{ineq-0209-02}).
We first rewrite ,
\begin{align*}
&\int_1^{\infty}\int_{|\bar{x}|>2}\left|\widetilde{K}_{k,t}(\bar{x}-\bar{y})-\widetilde{K}_{k,t}(\bar{x})\right|d\bar{x}\frac{dt}{t}\\
&=\int_1^{\infty}\int_{|\bar{x}|>2/t}\left|\widetilde{K}_{k}\left(\bar{x}-\frac{\bar{y}}{t}\right)-\widetilde{K}_{k}(\bar{x})\right|d\bar{x}\frac{dt}{t}\\
&\leq\int_1^{\infty}\omega_{1}(\widetilde{K}_{k},1/t)\frac{dt}{t}
=\int_0^{1}\omega_{1}(\widetilde{K}_{k},t)\frac{dt}{t},
\end{align*}
where $$\omega_{1}(\widetilde{K}_{k},t)=\sup_{|\bar{y}|=1}\int\left|\widetilde{K}_{k}\left(\bar{x}-t\bar{y}\right)-\widetilde{K}_{k}(\bar{x})\right|d\bar{x}.$$
Using (\ref{ineq-020904}), it can be deduced that
\begin{align*}
\omega_{1}(\widetilde{K}_{k},t)
&=\sup_{|\bar{y}|=1}\int\left|\int\sum_{i=1}^{n+1}2^{k(n+2)}(x_i-ty_i)\phi_i(2^k(x-ty-z,x_{n+1}-ty_{n+1}-(D+\Phi(z))))\right.\\
&\,\,\,\,\,\,\,\,\,\,\,\,\,\,\,\,\,\,\,-\left.\sum_{i=1}^{n+1}2^{k(n+2)}x_i\phi_i(2^k(x-z,x_{n+1}-(D+\Phi(z))))\rho(z)dz\right|d\bar{x}\\
&\leq\int\sup_{|\bar{y}|=1}\int\left|\sum_{i=1}^{n+1}2^{k(n+2)}(x_i-ty_i)\phi_i(2^k(x-ty-z,x_{n+1}-ty_{n+1}-(D+\Phi(z))))\right.\\
&\,\,\,\,\,\,\,\,\,\,\,\,\,\,\,\,\,\,\,-\left.\sum_{i=1}^{n+1}2^{k(n+2)}x_i\phi_i(2^k(x-z,x_{n+1}-(D+\Phi(z))))\right|d\bar{x}\rho(z)dz\\
&\leq \int\sup_{|\bar{y}|=1}\int\left|\sum_{i=1}^{n+1}2^{k(n+2)}(x_i-ty_i)\phi_{i}(2^k(\overline{x}-t\overline{y}))
-\sum_{i=1}^{n+1}2^{k(n+2)}x_i\phi_i(2^k\overline{x})
\right|d\bar{x}\rho(z)dz\\
&\,\,\,\,\,\,\,\,\,\,\,\,\,\,\,\,\,\,\,+\int\sup_{|\bar{y}|=1}\int\left|\sum_{i=1}^{n}2^{k(n+2)}z_i\phi_{i}(2^k(\overline{x}-t\overline{y}))\right.
+2^{k(n+2)}(D+\Phi(z))\phi_{n+1}(2^k(\overline{x}-t\overline{y}))\\
&\,\,\,\,\,\,\,\,\,\,\,\,\,\,\,\,\,\,\,-\left.\sum_{i=1}^{n}2^{k(n+2)}z_i\phi_i(2^k\overline{x})-2^{k(n+2)}(D+\Phi(z))\phi_{n+1}(2^k\overline{x})
\right|d\bar{x}\rho(z)dz\\
&\leq C(1+D)\max\{2^kt,\min\{2^{2k}t,2^{k}\}\}\leq C(1+D)\max\{1,2^k\}\min\{1,2^kt\}.
\end{align*}
Thus for $k\geq0$,
$$
\int_0^1\omega_{1}(\widetilde{K}_{k},t)\frac{dt}{t}\leq C(1+D)\max\{1,2^k\}\left(\int_{0}^{2^{-k}}2^kdt+\int_{2^{-k}}^{1}\frac{dt}{t}\right)\leq  C(1+D)2^k(1+k),
$$
and for $k<0$
$$
\int_0^1\omega_{1}(\widetilde{K}_{k},t)\frac{dt}{t}\leq C(1+D)\max\{1,2^k\}\int_{0}^{1}2^kdt\leq C(1+D)2^k.
$$
The proof is complete.

\subsection{The proof of Lemma \ref{le-0103}}
We first state a vector-valued variant of a theorem  about multipliers with localized bounds, which was introduced by Seeger.
Readers can find the  proof of  this  theorem in Jones-Seeger-Wright\cite{2008-strongvariation} and Guo\cite{2020-GUO-variation}.
\begin{prop}\label{pro-001}
Let $I\subset \mathbb{R}$ be a compact interval.
Let $\tilde{m}_{u}(\overline{\xi}):u\in I$ be a family of Fourier multipliers on
$\mathbb{R}^{n+1}$,
each of which is compactly supported on $\{\overline{\xi}:\frac12\leq |\overline{\xi}|\leq 2\}$
and satisfies
$$
\sup_{u\in I}\left|\partial^{\tau}_{\overline{\xi}}\tilde{m}_{u}(\overline{\xi})\right|\leq B\,\,\,\,\mbox{for each $0\leq |\tau|\leq n+2$}
$$
for some constant $B$.
For $u\in I$ and $j\in\mathbb{Z}$, denote by $T_{u,j}$
the multiplier operator with multiplier $\tilde{m}_{u}(2^{-j}\overline{\xi})$.
Fix some $p\in[2,\infty)$.
Assume that there exists some constant $A$ such that
$$
\sup_{j\in\mathbb{Z}}\left\|\left\|T_{u,j}f\right\|_{L^2(I)}\right\|_{L^s(\mathbb{R}^{n+1})}\leq
A\|f\|_{L^s(\mathbb{R}^{n+1})}
$$
for both $s=p$ and $s=2$.
Then
$$
\|\|\|T_{u,j}f\|_{L^2(I)}\|_{l^2(\mathbb{Z})}\|_{L^p(\mathbb{R}^{n+1})}\lesssim
A\left|\log\left(2+\frac{B}{A}\right)\right|^{\frac12-\frac1p}\|f\|_{L^p(\mathbb{R}^{n+1})}.
$$

\end{prop}

Firstly, we consider the case $k \leq0$.
Using inequality (\ref{ineq-sr-estimate}), we obtain
\begin{align*}
\|V_2^j(\mathcal{A}_kf)\|_{l^2_j}&\leq C \|\mathcal{A}_{k,t}f\|_{L^2_{t}}^{1/2}\left\|t\frac{d}{dt}\left(\mathcal{A}_{k,t}f\right)\right\|_{L^2_{t}}^{1/2}\\
&\leq C (1+D^{1/2})\|\mathcal{A}_{k,t}f\|_{L^2_{t}}+(1+D^{1/2})^{-1}\left\|t\frac{d}{dt}\left(\mathcal{A}_{k,t}f\right)\right\|_{L^2_{t}}.
\end{align*}
Recall that the kernel function of $\mathcal{A}_{k,t}$  is $H_{k,t}=K_{k,t}-G_{k,t}$.
As in the proof of Lemma \ref{le-0102},
we only prove
\begin{equation}\label{ineq-0910-02}
\|K_{k,t}*f\|_{L^p(L^2_{t})}\leq C2^{k\epsilon}\|f\|_{L^p}
\end{equation}
and
\begin{equation}\label{ineq-0910-03}
\left\|t\frac{d}{dt}K_{k,t}*f\right\|_{L^p(L^2_{t})}\leq C(1+D)2^{k\epsilon}\|f\|_{L^p}
\end{equation}
for $p>2$.

Using Minkowski's inequality, it can be deduced from
(\ref{ineq-0210-01}) that
$\|K_{k}\|_{L^1}\leq 1$.
Thus, by  H\"{o}lder inequality on $t$, we obtain
\begin{align*}
\left(\int_{\mathbb{R}^{n+1}}\left(\int_{1}^{2}|\mathcal{A}_{k,t}f(\overline{x})|^2dt\right)^{p/2}d\overline{x}\right)^{1/p}
&\leq \left(\int_{\mathbb{R}^{n+1}}\int_{1}^{2}|\mathcal{A}_{k,t}f(\overline{x})|^pdtd\overline{x}\right)^{1/p}\\
&=\left(\int_{1}^{2}\int_{\mathbb{R}^{n+1}}|\mathcal{A}_{k,t}f(\overline{x})|^pd\overline{x}dt\right)^{1/p}\\
&=\left(\int_{1}^{2}\int_{\mathbb{R}^{n+1}}|H_{k,t}*f(\overline{x}-tD)|^pd\overline{x}dt\right)^{1/p}\\
&\leq C\|f\|_{L^p}.
\end{align*}
Then, we get (\ref{ineq-0910-02})  if we  using Proposition \ref{pro-001} with $B=\mathcal{O}(2^{k(n+2)})$ and $A$ being the constant above.

Using Minkowski's inequality again, it can be deduced from
 (\ref{ineq-020903}) that
$\|\widetilde{K}_{k}\|_{L^1}\leq (1+D)$.
Note that $t\frac{d}{dt}K_{k,t}=-(n+1)K_{k,t}-\widetilde{K}_{k,t}$.
These results  imply (\ref{ineq-0910-03}).

It remains for us to consider the case $k\geq 1$.
Using inequality (\ref{ineq-sr-estimate}), we obtain
\begin{align*}
\|V_2^j(\mathcal{A}_kf)\|_{l^2_j}&\leq C \|\mathcal{A}_{k,t}f\|_{L^2_{t}}^{1/2}\left\|t\frac{d}{dt}\left(\mathcal{A}_{k,t}f\right)\right\|_{L^2_{t}}^{1/2}\\
&\leq C (1+D^{1/2})2^{\frac{k}{2}}\|\mathcal{A}_{k,t}f\|_{L^2_{t}}+(1+D^{1/2})^{-1}2^{-\frac{k}{2}}\left\|t\frac{d}{dt}\left(\mathcal{A}_{k,t}f\right)\right\|_{L^2_{t}}.
\end{align*}
Without loss of generality, using Proposition \ref{pro-001} again,  it remains for us to prove
\begin{equation}\label{ineq-0911-01}
\|\mathcal{A}_{k,t}f\|_{L^p\left(L^2_{t}([1,2])\right)}\leq C_{\epsilon}2^{-k\cdot\frac{n+1}{p}+k\epsilon}\|f\|_{L^p}
\end{equation}
and
\begin{equation}\label{ineq-0911-02}
\left\|t\frac{d}{dt}\mathcal{A}_{k,t}*f\right\|_{L^p\left(L^2_{t}([1,2])\right)}\leq C_{\epsilon}(1+D)2^{-k(\frac{n+1}{p}-1)+k\epsilon}\|f\|_{L^p}
\end{equation}
for $p>2+\frac{4}{n}$.

Now, let us focus our attention on the proof of (\ref{ineq-0911-01}).
Recall that the  multipliers associated  to operator $\mathcal{A}_{k,t}$
is $$m_{k}(\overline{\xi})=\mu(\overline{\xi})\psi(2^{-k}|\overline{\xi}|)-e^{-|\overline{\xi}|-iD\xi_{n+1}}\psi(2^{-k}|\overline{\xi}|),$$
and the expression of $\mu$ is
\begin{align*}
\mu(\overline{\xi})&=\int e^{-i\left(x\cdot\xi+\xi_{n+1}(D+\Phi(x))\right)}\rho(x)dx\\
&=e^{-iD\xi_{n+1}}\int e^{-i\xi_{n+1}\left(-x\cdot s+\Phi(x)\right)}\rho(x)dx
\end{align*}
where
$
s=\left(-\frac{\xi_1}{\xi_{n+1}},\ldots,-\frac{\xi_n}{\xi_{n+1}}\right)
$
for $\xi_{n+1}\neq0$.
Put $$\Phi^*(x,s)=-x\cdot s+\Phi(x).$$
Then we have
$$
\partial_{x_i}\Phi^*(x,s)=-s_i+\partial_{x_i}\Phi(x)
$$
and
$$
\partial_{x_ix_j}^2\Phi^*(x,s)=\partial_{x_ix_j}^2\Phi(x).
$$

Recall that $\Phi(x)=\sum_{i=1}^{n}\lambda_ix_i^2+R(x)$.
Thus we can get
\begin{equation}\label{equ021301}
\partial_{x_ix_j}^2\Phi^*(0,s)=\left\{\begin{array}{ll}
2\lambda_i,& i=j\,,\\
0,& i\neq j\,.
\end{array}\right.
\end{equation}
Since $|x|$ is assumed to be sufficiently small,
then the implicit function theorem and (\ref{equ021301}) imply  that
there exists a smooth solution
$x^*=x^*(s)$ of the equation  $\partial_{x_i}\Phi^*(x,s)=0$ for $i=1,\ldots,n$.
From the above arguments and Taylor expansion of the smooth functions $\Phi(x)$,
the function $x^*(s)$ can be written as
\begin{equation}\label{equ-x-ex}
x^*(s)=\left(\frac{s_1}{2\lambda_1},\ldots,
\frac{s_n}{2\lambda_n}
\right)+R_{ph}(s)
\end{equation}
where
$R_{ph}(0)=\left.D_{s}R_{ph}(s)\right|_{s=0}=0$.

Denote
$\Phi^{**}(s)=-s\cdot x^*(s)+\Phi(x^*(s))$.
By applying the method of stationary phase, we obtain
\begin{equation}
\mu(\overline{\xi})=e^{-iD\xi_{n+1}}\left(e^{i\xi_{n+1}\Phi^{**}(s)}\frac{A_{m}(\overline{\xi})}{(1+|\overline{\xi}|)^{\frac{n}{2}}}\chi(s)
+A_{re}(\overline{\xi})\right).
\end{equation}
where $\chi(s)$ is a smooth function supported in $[0,c]$, $A_{m}(\overline{\xi})$ is a symbol of order of zero.
More precisely, $A_{m}(\overline{\xi})$ satisfies
$$
\left|D^{\overline{\alpha}}_{\overline{\xi}}A_m(\overline{\xi})\right|\leq C\left(1+|\overline{\xi}|\right)^{-|\overline{\alpha}|}.
$$
Moreover, the remainder term $A_{re}(\overline{\xi})$ satisfies
$$
\left|D^{\overline{\alpha}}_{\overline{\xi}}A_{re}(\overline{\xi})\right|\leq C_{N,\overline{\alpha}}\left(1+|\overline{\xi}|\right)^{-N}.
$$
At last, the multiplier corresponding to $\mathcal{A}_{k,t}$ can be written by
\begin{equation*}
m_k(\overline{\xi})
=e^{-iD\xi_{n+1}}\left(e^{i\xi_{n+1}\Phi^{**}(s)}\frac{A_{m}(\overline{\xi})}{(1+|\overline{\xi}|)^{\frac{n}{2}}}\chi(s)\right)\psi(2^{-k}\overline{\xi})
+e^{-iD\xi_{n+1}}A_{re}(\overline{\xi})\psi(2^{-k}\overline{\xi}),
\end{equation*}
where
$$
\left|D^{\overline{\alpha}}_{\overline{\xi}}A_{re}(\overline{\xi})\right|\leq C_{N,\overline{\alpha}}\left(1+|\overline{\xi}|\right)^{-N}.
$$

Let
$$
\mathcal{M}_{k,t}f(\overline{x})=\frac{1}{(2\pi)^{n+1}}\int_{\mathbb{R}^{n+1}}e^{it(\overline{x}\cdot\overline{\xi}-D\xi_{n+1}+\xi_{n+1}\Phi^{**}(s))}
\frac{A_{m}(\overline{\xi})}{(1+|\overline{\xi}|)^{\frac{n}{2}}}\chi(s)\psi(2^{-k}\overline{\xi})\widehat{f}(\overline{\xi})d\overline{\xi}
$$
denote the main term and let
$$
\mathcal{B}_{k,t}f(\overline{x})=\frac{1}{(2\pi)^{n+1}}\int_{\mathbb{R}^{n+1}}e^{it(\overline{x}\cdot\overline{\xi}-D\xi_{n+1})}
A_{re}(\overline{\xi})\psi(2^{-k}\overline{\xi})\widehat{f}(\overline{\xi})d\overline{\xi}
$$
denote the remainder term.
This reduces the matter to proving the associated desired estimates for $\mathcal{M}_{k,t}$ and $\mathcal{B}_{k,t}$.

\textit{The estimate of the remainder term.}
Denote the kernel function of $\mathcal{B}_{k,t}$ by
 \begin{equation*}
 K_{ek}(\overline{x},t)=\frac{1}{(2\pi)^{n+1}}\int_{\mathbb{R}^{n+1}}e^{it(\overline{x}\cdot\overline{\xi}-D\xi_{n+1})}
A_{re}(\overline{\xi})\psi(2^{-k}\overline{\xi})d\overline{\xi}.
 \end{equation*}
Integrating by parts in $\overline{\xi}$, since $t\approx 1$, we obtain
\begin{eqnarray*}
\left|K_{ek}(\overline{x},t)\right|\lesssim_N2^{-kN}\frac{2^{k(n+1)}}{1+|2^kx|^N+(2^k|x_{n+1}-D|)^N}
\end{eqnarray*}
Using H\"{o}lder's inequality and Young's inequality for convolution, it can be deduced that
\begin{equation}
\left\|\mathcal{B}_{k,t}f\right\|_{L^p(L^2([1,2]))}=\left\|K_{ek}(\overline{x},t)*_{\overline{x}}f\right\|_{L^p(L^2([1,2]))}\lesssim_{N} 2^{-kN}\|f\|_{L^p}
\end{equation}
for $p\geq2$.
In order to prove (\ref{ineq-0911-01}),
it remains for us to consider the estimate associated to the  main term.

\textit{The estimate of the main term.}
Before we begin our proof, we introduce two inequalities about local smoothing estimates which can be found in the proof of Theorem 3.1
and
Theorem 3.2 in \cite{1993-MSS-loacal-smoothing}.
\begin{lemma}\label{le-010301}
Suppose $a(\overline{\xi},\overline{x})$ is a symbol of order zero.
That is
$$
\left|\left(\frac{\partial}{\partial\overline{\xi}}\right)^{\alpha}
\left(\frac{\partial}{\partial \overline{x}}\right)^{\beta}
a(\overline{\xi},\overline{x})\right|\leq C_{\alpha,\beta}(1+|\overline{\xi}|)^{-|\alpha|}.
$$
Denote
\begin{equation}\label{ineq-local}
\mathcal{Q}_{\lambda}f(\overline{x},t)=\widetilde{\psi}(t)\int_{\mathbb{R}^{n+1}}e^{it(\overline{x}\cdot\overline{\xi}-h(\overline{\xi}))}
a_{\lambda}(\overline{\xi},\overline{x},\lambda)\widehat{f}(\overline{\xi})d\overline{\xi},
\end{equation}
where
$a_{\lambda}(\overline{\xi},\overline{x})=a(\overline{\xi},\overline{x})\psi\left(\lambda^{-1}|\overline{\xi}|\right)$
,
$h(\cdot)$ is homogeneous of degree 1
and
$\widetilde{\psi}(t)$ is a positive smooth function supported in $[1/2,5/2]$ and
$\widetilde{\psi}(t)=1$ for $t\in[1,2]$.
Then we have
\begin{equation}
\left\|\mathcal{Q}_{\lambda}f\right\|_{L^p\left(\mathbb{R}^{n+1},L^2(\mathbb{R})\right)}\leq \lambda^{(n+1)\left(\frac12-\frac1p\right)-\frac12}
\|f\|_{L^2(\mathbb{R}^{n+1})},
\end{equation}
and
\begin{equation}
\left\|\mathcal{Q}_{\lambda}f\right\|_{L^p\left(\mathbb{R}^{n+2})\right)}\leq \lambda^{(n+1)\left(\frac12-\frac1p\right)-\frac1p}
\|f\|_{L^2(\mathbb{R}^{n+1})},
\end{equation}
for $p>2+\frac{4}{n}$.
\end{lemma}

Denote
$$
\mathcal{M}^*_{k,t}f(\overline{x})=\frac{1}{(2\pi)^{n+1}}\int_{\mathbb{R}^{n+1}}e^{it(\overline{x}\cdot\overline{\xi}+\xi_{n+1}\Phi^{**}(s))}
\frac{A_{m}(\overline{\xi})}{(1+|\overline{\xi}|)^{\frac{n}{2}}}\chi(s)\psi(2^{-k}|\overline{\xi}|)\widehat{f}(\overline{\xi})d\overline{\xi}.
$$
In order to use Lemma \ref{le-010301}, we let $h(\overline{\xi})=-\xi_{n+1}\Phi^{**}(s)$ and
$$ a(\overline{\xi},\overline{x}) =2^{k\cdot\frac{n}{2}}\frac{A_{m}(\overline{\xi})}{(1+|\overline{\xi}|)^{\frac{n}{2}}}\chi(s)\psi(2^{-k}|\overline{\xi}|)$$
since $|\overline{\xi}|\approx 2^k$.
It follows that,
\begin{equation}\label{ineq-main-A-001}
\left\|\mathcal{M}_{k,t}^*f\right\|_{L^p\left(\mathbb{R}^{n+1},L^2([1,2])\right)}\leq C 2^{-k\cdot\frac{n+1}{p}}
\|f\|_{L^2(\mathbb{R}^{n+1})},
\end{equation}
and
\begin{equation}\label{ineq-main-A-0214}
\left\|\mathcal{M}_{k,t}^*f\right\|_{L^p\left(\mathbb{R}^{n+1},L^p([1,2])\right)}\leq C
2^{k\left((n+1)\left(\frac12-\frac1p\right)-\frac1p-\frac{n}{2}\right)}
\|f\|_{L^2(\mathbb{R}^{n+1})},
\end{equation}
for $p>2+\frac{4}{n}$.

Recalling that $h(\cdot)$ is homogeneous of degree 1,
we have $h(\overline{\xi})=\sum_{i=1}^{n+1}\xi_i(\partial_ih)\left(\frac{\overline{\xi}}{|\overline{\xi}|}\right)$.
Thus we can rewrite $\mathcal{M}^*_{k,t}f$ as a convolution, that is,
\begin{equation}\label{ineq-010501}
\mathcal{M}^*_{k,t}f(\overline{x})=(f*_{\overline{x}}K^*_{m,k})(\overline{x},t),
\end{equation}
where the kernel function $K^*_{m,k}$ is given by
\begin{equation}
K^*_{m,k}(\overline{x},t)=\frac{1}{(2\pi)^{n+1}}\int_{\mathbb{R}^{n+1}}e^{it\left(\sum_{i=1}^{n+1}\xi_i\left(x_i-(\partial_ih)
\left(\frac{\overline{\xi}}{|\overline{\xi}|}\right)\right)\right)}
\frac{A_{m}(\overline{\xi})}{(1+|\overline{\xi}|)^{\frac{n}{2}}}\chi(s)\psi(2^{-k}|\overline{\xi}|)d\overline{\xi}.
%&=&\frac{2^{k(n+1)}}{(2\pi)^{n+1}}\int_{\mathbb{R}^{n+1}}e^{it2^k\left(\sum_{i=1}^{n+1}\xi_i\left(x_i-(\partial_ih)
%\left(\frac{\overline{\xi}}{|\overline{\xi}|}\right)\right)\right)}
%\frac{A_{m}(2^k\overline{\xi})}{(1+2^k|\overline{\xi}|)^{\frac{n}{2}}}\chi(s)\psi(|\overline{\xi}|)d\overline{\xi}\nonumber\\
\end{equation}
Let $\mathcal{C}=\sup\{|(\partial_i h)(\overline{\theta})|:{\overline{\theta}\in\mathbb{S}^n},i=1,\ldots,n+1\}$.
By integrating by parts in $\overline{\xi}$ and $t\approx1$, we can infer that for any $N\in\mathbb{N}$,
\begin{equation}\label{ineq-010502}
|K^*_{m,k}(\overline{x},t)|\lesssim_N 2^{-kN}  \frac{2^{k(n+1)}}{1+|2^k\overline{x}|^N}
\end{equation}
whenever $|\overline{x}|\geq 2(n+1)\mathcal{C}$.
%{\color{red} \remark It is easy to show for $|x|>2a$, we have
%$$
%\frac{1}{|\lambda x|^{2N}}\leq \frac{1}{(2a\lambda)^{N}}\frac{1}{|\lambda x|^{N}}.
%$$
%Then we can get the result by the following inequality
%$$
%\frac{1}{|\lambda x|^{N}}\leq \frac{1+\frac{1}{(2a\lambda)^N}}{1+|\lambda x|^N}.
%$$
%}

For $\overline{z}:=(z_1,\cdots,z_{n+1})$,
denote by $Q_{\overline{z}}$ the cube $\prod_{i=1}^{n+1}[z_i,z_i+\mathcal{C})$
and denote by $Q_{\overline{z}^*}$ the cube centered at $\overline{z}$ with side-length $2(n+1)\mathcal{C}$.
By using the notation $f_{Q_{\overline{z}^*}}:=f1_{Q_{\overline{z}^*}}$,
and recalling the expression (\ref{ineq-010501}),
we can rewrite
$\mathcal{M}^*_{k,t}f(\overline{x})$ as
\begin{eqnarray}
\sum_{\overline{z}}(f_{Q_{\overline{z}}}*_{\overline{x}}K^*_{m,k})(\overline{x},t)
&=&\sum_{\overline{z}}(f_{Q_{\overline{z}}}*_{\overline{x}}K^*_{m,k})(\overline{x},t)1_{Q_{\overline{z}}^*}(\overline{x})\nonumber\\
&&\,\,\,\,\,\,\,+\sum_{\overline{z}}(f_{Q_{\overline{z}}}*_{\overline{x}}K^*_{m,k})(\overline{x},t)1_{(Q_{\overline{z}}^*)^c}(\overline{x})\nonumber\\
&=:&\mathcal{T}_{1mk}(\overline{x},t)+\mathcal{T}_{2mk}(\overline{x},t).\nonumber
\end{eqnarray}

By (\ref{ineq-main-A-001}), we can get the estimate
\begin{eqnarray}\label{ineq-021401}
\|\mathcal{T}_{1mk}\|_{L^p(L^2_t([1,2]))} &\lesssim&
\left(\sum_{\overline{z}}\left\|\mathcal{M}^*_{k,t}f_{Q_{\overline{z}}}\right\|_{L^p(L^2_t([1,2]))}^p\right)^{1/p}\nonumber\\
&\lesssim&2^{-k\cdot\frac{n+1}{p}}\left(\sum_{\overline{z}}\|f_{Q_{\overline{z}}}\|_{L^2}^p\right)^{1/p}.
\end{eqnarray}
Using Fubini's theorem, we obtain $\|f_{Q_{\overline{z}}}\|_{L^2}\lesssim \|f_{Q_{\overline{z}}}\|_{L^p}$.
It follows from (\ref{ineq-021401}) that
\begin{equation}\label{ineq-021402}
\|\mathcal{T}_{1mk}\|_{L^p(L^2_t([1,2]))} \lesssim
2^{-k\cdot\frac{n+1}{p}}\|f_{Q_{\overline{z}}}\|_{L^p}.
\end{equation}

On the other hand,
using (\ref{ineq-010502}), we immediately deduce the pointwise inequality
\begin{equation}\label{ineq-01050201}
|\mathcal{T}_{2mk}|(\overline{x},t)\lesssim_N 2^{-kN}|f|*\mathcal{E}_{k,N,t}(\overline{x}),
\end{equation}
where
$$
\mathcal{E}_{k,N,t}=\frac{2^{k(n+1)}t^{n+1}}{1+|2^kt\overline{x}|^N}.
$$
Then, taking the $L^p(L^2_t([1,2]))$ norm on both sides of (\ref{ineq-01050201}),
and using H\"{o}lder's inequality and Young's inequality for convolution, we can obtain
\begin{equation}\label{ineq-021403}
\left\|\mathcal{T}_{2mk}\right\|_{L^p(L^2_t([1,2]))}
\lesssim_N 2^{-kN}\|f\|_{L^p}.
\end{equation}
It can be deduced from (\ref{ineq-021402}) and (\ref{ineq-021403}) that
\begin{equation}\label{ineq-0214-05}
\left\|\mathcal{M}_{k,t}^*f\right\|_{L^p\left(\mathbb{R}^{n+1},L^2([1,2])\right)}\leq C_{\epsilon} 2^{-k\cdot\frac{n+1}{p}+k\epsilon}
\|f\|_{L^p(\mathbb{R}^{n+1})}.
\end{equation}
After going through the same steps to prove (\ref{ineq-0214-05}), we obtain
\begin{equation}\label{ineq-0214-06}
\left\|\mathcal{M}_{k,t}^*f\right\|_{L^p\left(\mathbb{R}^{n+1},L^p([1,2])\right)}\leq C_{\epsilon}
2^{k\left((n+1)\left(\frac12-\frac1p\right)-\frac1p+\frac{n}{2}+\epsilon\right)}
\|f\|_{L^p(\mathbb{R}^{n+1})}.
\end{equation}
Note that $\mathcal{M}_{k,t}f(\overline{x})=\mathcal{M}^*_{k,t}f(x,x_{n+1}-D)$.
Thus we can get from (\ref{ineq-0214-05}),
\begin{equation}\label{ineq-0214-0501}
\left\|\mathcal{M}_{k,t}f\right\|_{L^p\left(\mathbb{R}^{n+1},L^2([1,2])\right)}\leq C_{\epsilon,D} 2^{-k\cdot\frac{n+1}{p}+k\epsilon}
\|f\|_{L^p(\mathbb{R}^{n+1})}.
\end{equation}
On the other hand, using H\"{o}lder's inequality, we have
\begin{equation*}
\left\|\mathcal{M}_{k,t}f\right\|_{L^p\left(\mathbb{R}^{n+1},L^2([1,2])\right)}\leq\left\|\mathcal{M}_{k,t}f\right\|_{L^p\left(\mathbb{R}^{n+1},L^p([1,2])\right)}.
\end{equation*}
Combining the inequality above with (\ref{ineq-0214-06}), it can be deduced that
\begin{equation}\label{ineq-0214-06-01}
\left\|\mathcal{M}_{k,t}f\right\|_{L^p\left(\mathbb{R}^{n+1},L^2([1,2])\right)}\leq C_{\epsilon}
2^{k\left((n+1)\left(\frac12-\frac1p\right)-\frac1p+\frac{n}{2}+\epsilon\right)}
\|f\|_{L^p(\mathbb{R}^{n+1})}.
\end{equation}
It implies that the constant $C_{D,\epsilon}$ is bounded by a constant $C_{\epsilon}$ which is independent of $D$.
Then we obtained (\ref{ineq-0911-01}).

\textit{The proof of (\ref{ineq-0911-02}).}
Denote the multiplier associated to $t\frac{d}{dt}\left(\mathcal{A}_{k,t}f\right)$ by $\tilde{m}_k(t\overline{\xi})$.
Then
\begin{equation}\label{ineq-t-001}
\tilde{m}_k(\overline{\xi})=\sum_{i=1}^{n+1}\xi_{i}(\partial_{\xi_i}m)(\overline{\xi})\psi(2^{-k}|\overline{\xi}|)
+2^{-k}\sum_{i=1}^{n+1}\frac{\xi_{i}^2}{|\overline{\xi}|}m(\overline{\xi})(\partial_{\xi_i}\psi)(2^{-k}|\overline{\xi}|).
\end{equation}
Note that
\begin{equation}\label{ineq-t-002}
(\partial_{\xi_i}\mu)(\overline{\xi})=-i\int e^{-i\left(x\cdot\xi+\xi_{n+1}(D+\Phi(x))\right)}x_i\rho(x)dx
\,,\,\,\,\,\,\,\,\mbox{for $i=1,\ldots,n$}
\end{equation}
and
\begin{equation}\label{ineq-t-003}
(\partial_{\xi_{n+1}}\mu)(\overline{\xi})=-i\int e^{-i\left(x\cdot\xi+\xi_{n+1}(D+\Phi(x))\right)}(D+\Phi(x))\rho(x)dx.
\end{equation}
Recall that
$m(\overline{\xi})=\mu(\overline{\xi})-e^{-iD\xi_{n+1}-|\overline{\xi}|}$.
By applying the method of stationary phase, it follows by (\ref{ineq-t-001}), (\ref{ineq-t-002})
and (\ref{ineq-t-003}) that
\begin{equation}
\tilde{m}_k(\overline{\xi})
=(1+D)2^ke^{-iD\xi_{n+1}}\left(e^{i\xi_{n+1}\Phi^{**}(s)}\frac{\widetilde{A}_{m}(\overline{\xi})}{(1+|\overline{\xi}|)^{\frac{n}{2}}}\tilde{\chi}(s)\right)\psi(2^{-k}\overline{\xi})
+e^{-iD\xi_{n+1}}\widetilde{A_{re}}(\overline{\xi})\psi(2^{-k}\overline{\xi}),
\end{equation}
where $\tilde{\chi}(s)$ is a smooth function supported in $[0,\tilde{c}]$, $\widetilde{A}_{m}(\overline{\xi})$ is a symbol of order of zero.
Moreover, the remainder term $\widetilde{A_{re}}(\overline{\xi})$ satisfies
$$
\left|D^{\overline{\alpha}}_{\overline{\xi}}\widetilde{A_{re}}(\overline{\xi})\right|\leq C_{N,\overline{\alpha}}\left(1+|\overline{\xi}|\right)^{-N}.
$$
Using the same method as in proving (\ref{ineq-0911-01}), we can obtain (\ref{ineq-0911-02}).
The proof is complete.

\subsection{The proof of Lemma \ref{le-0104}.}
Denote $\delta_{\rho}f(\overline{x})=f\left(\rho\overline{x}\right)$.
Since
\begin{align*}
\mathcal{A}_{k,2^jt}f(\overline{x})&=\frac{1}{(2\pi)^2}\int e^{i\overline{x}\cdot\overline{\xi}}m(2^{j}t\overline{\xi})
\psi(2^{-k+j}t|\overline{\xi}|)\widetilde{\psi}(2^{-k+j}|\overline{\xi}|)\widehat{f}(\overline{\xi})d\overline{\xi}\\
&=\frac{1}{(2\pi)^2}\int e^{i2^{-j}\overline{x}\cdot\overline{\xi}}m(t\overline{\xi})
\psi(2^{-k}t|\overline{\xi}|)\mathcal{F}\left(\delta_{2^j}(\tilde{P}_{k-j}f)\right)(\overline{\xi})d\overline{\xi}\\
&=\mathcal{A}_{k,t}\left(\delta_{2^j}(\tilde{P}_{k-j}f)\right)(2^{-j}|\overline{x}|),
\end{align*}
we have
$V_4^j(\mathcal{A}_kf)(\overline{x})=V_4^1\left(\mathcal{A}_k\left(\delta_{2^j}(\tilde{P}_{k-j}f)\right)\right)(2^{-j}\overline{x})$.
Thus if we prove
\begin{equation}\label{ineq-short-44-021401}
\left\|V_4^{1}(\mathcal{A}_kf)\right\|_{L^4}\leq C_{\epsilon}(1+D^{1/4}) \min\left\{2^{-k\left(\frac{1}{8}-\epsilon\right)},1\right\}\|f\|_{L^4},
\end{equation}
then we obtain
\begin{align*}
\left\|V_4^{j}(\mathcal{A}_kf)\right\|_{L^4(l^4_j)}
&=\left\|V_4^1\left(\mathcal{A}_k\left(\delta_{2^j}(\tilde{P}_{k-j}f)\right)\right)(2^{-j}\cdot)\right\|_{L^4(l^4_j)}\\
&\leq C_{\epsilon}(1+D^{1/4}) \min\left\{2^{-k\left(\frac{1}{8}-\epsilon\right)},1\right\}
\left\|(\tilde{P}_{k-j}f)\right\|_{L^4(l^4_j)}\\
&\leq C_{\epsilon}(1+D^{1/4}) \min\left\{2^{-k\left(\frac{1}{8}-\epsilon\right)},1\right\}\|f\|_{L^4}.
\end{align*}
It remains for us to prove (\ref{ineq-short-44-021401}).

Using inequality (\ref{ineq-sr-estimate}), we obtain
\begin{equation*}
V_4^j(\mathcal{A}_kf)\leq C \|\mathcal{A}_{k,t}f\|_{L^4_{t}([1,2])}^{3/4}\left\|t\frac{d}{dt}\left(\mathcal{A}_{k,t}f\right)\right\|_{L^4_{t}([1,2])}^{1/4}.
\end{equation*}
By H\"{o}lder's inequality, we obtain
\begin{equation*}
\|V_4^j(\mathcal{A}_kf)\|_{L^4}\leq C \|\mathcal{A}_{k,t}f\|_{L^4(L^4_{t}([1,2]))}^{3/4}\left\|t\frac{d}{dt}\left(\mathcal{A}_{k,t}f\right)\right\|_{L^4(L^4_{t}([1,2]))}^{1/4}.
\end{equation*}
It is enough for us to prove
\begin{equation}\label{ineq-0214-11}
\|\mathcal{A}_{k,t}f\|_{L^4\left(L^4_{t}([1,2])\right)}\leq C_{\epsilon}\min\left\{2^{-\frac{3}{8}\cdot  k+k\epsilon},1\right\}\|f\|_{L^4}
\end{equation}
and
\begin{equation}\label{ineq-0214-12}
\left\|t\frac{d}{dt}\mathcal{A}_{k,t}*f\right\|_{L^4\left(L^4_{t}([1,2])\right)}\leq C_{\epsilon}(1+D)
\left\{2^{\frac{5}{8}\cdot k+k\epsilon},1\right\}\|f\|_{L^4}.
\end{equation}

Recall that
the multipliers corresponding to $\mathcal{A}_{k,t}$ and $t\frac{d}{dt}\mathcal{A}_{k,t}$  can be written as
\begin{equation*}
m_k(\overline{\xi})
=e^{-iD\xi_{n+1}}\left(e^{i\xi_{n+1}\Phi^{**}(s)}\frac{A_{m}(\overline{\xi})}{(1+|\overline{\xi}|)^{\frac{n}{2}}}\chi(s)\right)\psi(2^{-k}\overline{\xi})
+e^{-iD\xi_{n+1}}A_{re}(\overline{\xi})\psi(2^{-k}\overline{\xi}),
\end{equation*}
and \begin{equation*}
\tilde{m}_k(\overline{\xi})=\sum_{i=1}^{n+1}\xi_{i}(\partial_{\xi_i}m)(\overline{\xi})\psi(2^{-k}|\overline{\xi}|)
+2^{-k}\sum_{i=1}^{n+1}\frac{\xi_{i}^2}{|\overline{\xi}|}m(\overline{\xi})(\partial_{\xi_i}\psi)(2^{-k}|\overline{\xi}|)
\end{equation*}
respectively.
Using the same method as in the proof of  Lemma \ref{le-0103}, we can get
(\ref{ineq-0214-11})  and (\ref{ineq-0214-12}) by the following lemma,  introduced in the proof of Theorem 6.1 in \cite{1993-MSS-loacal-smoothing}.
\begin{lemma}
Denote $\mathcal{Q}_{\lambda}f(\overline{x},t)$ as in Lemma \ref{le-010301}.
Then we have
\begin{equation}
\left\|\mathcal{Q}_{\lambda}f\right\|_{L^4\left(\mathbb{R}^{3}\right)}\leq C_{\epsilon} \lambda^{\frac18+\epsilon}
\|f\|_{L^4(\mathbb{R}^{2})}.
\end{equation}
\end{lemma}

\subsection{The proof of the necessary condition}
In this subsection,
we establish the necessity of the condition $p<r(n+1)$ in Theorem \ref{theorem-k-principal} for the case $k=n$.
The argument follows a standard Knapp-type counterexample construction, adapted to the variational setting.

\textit{Preliminary work.}
Recall that the Fourier transform of the surface measure is given by
$$
\mu(\xi,\xi_{n+1})=e^{-iD\xi_{n+1}}\left(e^{i\xi_{n+1}\Phi^{**}(s)}\frac{A_{m}(\overline{\xi})}{(1+|\overline{\xi}|)^{\frac{n}{2}}}\chi(s)+A_{re}(\overline{\xi})\right).
$$
We split the operator $A_t$ into a main part and a remainder part,
\[
A_tf = Ma_tf + Er_tf,
\]
where
$$
Ma_tf(\overline{x})=\int_{\mathbb{R}^{n+1}}e^{-itD\xi_{n+1}}e^{it\xi_{n+1}\Phi^{**}(s)}\frac{A_{m}(\overline{\xi})}{(1+|\overline{\xi}|)^{\frac{n}{2}}}\chi(s)
\widehat{f}(\xi,\xi_{n+1})e^{i(x\cdot\xi+x_{n+1}\xi_{n+1})}d\overline{\xi},
$$
and
$$
Er_tf(\overline{x})=\int_{\mathbb{R}^{n+1}}e^{-itD\xi_{n+1}}A_{re}(\overline{\xi})
\widehat{f}(\xi,\xi_{n+1})e^{i(x\cdot\xi+x_{n+1}\xi_{n+1})}d\overline{\xi}.
$$
Consequently,
\begin{equation}\label{ineq-ne-001}
\left\|V_r\left(A_tf:t>0\right)\right\|_{L^p}\geq\left\|V_r\left(Ma_tf:t>0\right)\right\|_{L^p}-
\left\|V_r\left(Er_tf:t>0\right)\right\|_{L^p}.
\end{equation}
Thus it suffices to obtain a lower bound for the main term and an upper bound for the remainder term, which will be negligible for large $\lambda$.

To proceed, we perform a change of variables that linearises the phase in $t$. Set
\[
\xi_i=-\frac{s_i s_{n+1}}{-D+\Phi^{**}(s)}\qquad (i=1,\ldots,n),\qquad
\xi_{n+1}=\frac{s_{n+1}}{-D+\Phi^{**}(s)}.
\]
Then the Jacobian of this transformation is $|s_{n+1}|^n\Psi^*(s)$
for some smooth function $\Psi^*$ arising from the change of variables.
After this substitution, the main term and remainder term become
\begin{align*}
Ma_tf(\overline{x})&=\int_{\mathbb{R}^{n+1}}e^{its_{n+1}}\frac{A_{m}\left(-\frac{s_{n+1}\cdot s}{-D+\Phi^{**}(s)},\frac{s_{n+1}}{-D+\Phi^{**}(s)}\right)}{\bigl(1+\bigl|\bigl(-\frac{s_{n+1}\cdot s}{-D+\Phi^{**}(s)},\frac{s_{n+1}}{-D+\Phi^{**}(s)}\bigr)\bigr|\bigr)^{\frac{n}{2}}}\chi(s)\\
&\qquad\times e^{is_{n+1}\frac{-s\cdot x+x_{n+1}}{-D+\Phi^{**}(s)}}\widehat{f}\!\left(-\frac{s_{n+1}\cdot s}{-D+\Phi^{**}(s)},\frac{s_{n+1}}{-D+\Phi^{**}(s)}\right)|s_{n+1}|^n\Psi^*(s)\,d\overline{s},
\end{align*}
and similarly
\begin{align*}
Er_tf(\overline{x})&=\int_{\mathbb{R}^{n+1}}e^{its_{n+1}}A_{re}\!\left(-\frac{s_{n+1}\cdot s}{-D+\Phi^{**}(s)},\frac{s_{n+1}}{-D+\Phi^{**}(s)}\right)\\
&\qquad\times e^{is_{n+1}\frac{-s\cdot x+x_{n+1}}{-D+\Phi^{**}(s)}}\widehat{f}\!\left(-\frac{s_{n+1}\cdot s}{-D+\Phi^{**}(s)},\frac{s_{n+1}}{-D+\Phi^{**}(s)}\right)|s_{n+1}|^n\Psi^*(s)\,d\overline{s}.
\end{align*}

Now we choose a family of test functions $f_\lambda$ concentrated in frequency near a point where the phase is stationary.
Define $\widehat{f_\lambda}$ by
 \begin{equation}\label{test-function}
 \widehat{f_{\lambda}}\left(-\frac{s_{n+1}\cdot s}{-D+\Phi^{**}(s)},\frac{s_{n+1}}{-D+\Phi^{**}(s)}\right)=\chi_{0}(s)\chi_1\left(\frac{s_{n+1}}{\lambda}\right)e^{-i\frac{s_{n+1}^2}{2\lambda}}.
 \end{equation}
 Here $\chi_0$ is a smooth cut-off supported in a small neighbourhood of $s=0$,
 and $\chi_1$ is a smooth function supported in $[b/2,2b]$ for a small positive constant $c$.
 This choice makes $f_\lambda$ essentially a wave packet localized in both space and frequency; the quadratic exponential factor will later create oscillations that produce a large variation.

We first dispose of the remainder term. Substituting the definition of $\widehat{f_\lambda}$ into $Er_tf_\lambda$ gives
\begin{align*}
Er_tf_{\lambda}(\overline{x})&=\lambda^{n+1}\int_{\mathbb{R}^{n+1}}e^{it\lambda s_{n+1}}A_{re}\!\left(-\frac{\lambda s_{n+1}\cdot s}{-D+\Phi^{**}(s)},
\frac{\lambda s_{n+1}}{-D+\Phi^{**}(s)}\right)\\
&\qquad\times e^{i\lambda s_{n+1}\frac{-s\cdot x+x_{n+1}}{-D+\Phi^{**}(s)}}\chi_{0}(s)\chi_1(s_{n+1})e^{-i\frac{\lambda s_{n+1}^2}{2}}|s_{n+1}|^n\Psi^*(s)\,d\overline{s}.
\end{align*}
Because $A_{re}$ satisfies the rapid decay estimate
\begin{equation}\label{est-021601}
\left|D^{\overline{\alpha}}_{\overline{\xi}}A_{re}(\overline{\xi})\right|\leq C_{N,\overline{\alpha}}\left(1+|\overline{\xi}|\right)^{-N},
\end{equation}
a routine integration by parts argument (using the factor $e^{it\lambda s_{n+1}}$) shows that $Er_tf_\lambda$ is negligible. More precisely, for any $r>2$ and $p>1$,
\begin{equation}\label{usfunn-001-0}
\|V_r(Er_tf_\lambda:t>0)\|_{L^p}\lesssim_N \lambda^{-N}\|f_\lambda\|_{L^p}.
\end{equation}

Now we turn to the main term and estimate $\|f_\lambda\|_{L^p}$ and the variational norm of $Ma_tf_\lambda$.

\noindent\textit{Proof of (\ref{ineq-f-la}): upper bound for $\|f_\lambda\|_{L^p}$.}
By Fourier inversion,
\[
f_\lambda(\overline{x})=\frac{1}{(2\pi)^{n+1}}\int e^{i\overline{\xi}\cdot\overline{x}}\widehat{f_\lambda}(\overline{\xi})\,d\overline{\xi}.
\]
Changing variables as before,
\[
f_\lambda(\overline{x})=\frac{1}{(2\pi)^{n+1}}\int
e^{is_{n+1}\frac{-s\cdot x+x_{n+1}}{-D+\Phi^{**}(s)}}
\widehat{f_{\lambda}}\!\left(-\frac{s_{n+1}\cdot s}{-D+\Phi^{**}(s)},\frac{s_{n+1}}{-D+\Phi^{**}(s)}\right)
|s_{n+1}|^n\Psi^*(s)\,d\overline{s}.
\]
Inserting the definition (\ref{test-function}) yields
\begin{align}
f_{\lambda}(\overline{x})&=\frac{\lambda^{n+1}}{(2\pi)^{n+1}}\int
e^{i\lambda\left( s_{n+1}\frac{-s\cdot x+x_{n+1}}{-D+\Phi^{**}(s)}-\frac{s_{n+1}^2}{2}\right)}
\chi_{0}(s)\chi_1(s_{n+1})
|s_{n+1}|^n\Psi^*(s)\,d\overline{s}\nonumber\\
&=\frac{\lambda^{n+1}}{(2\pi)^{n+1}}\int
e^{i\lambda\left(-\frac{s_{n+1}^2}{2}\right)}J\bigl(\lambda x_{n+1}s_{n+1},\dot{x}\bigr)
\chi_1(s_{n+1})\,|s_{n+1}|^n\,ds_{n+1},
\end{align}
where $\dot{x}=(\frac{x_1}{x_{n+1}},\ldots,\frac{x_n}{x_{n+1}})$ and
\begin{equation}
J(\mu,\dot{x})=\int e^{i\mu\frac{-s\cdot\dot{x}+1}{-D+\Phi^{**}(s)}}\chi_0(s)\Psi^*(s)\,ds.
\end{equation}
The phase function $\Lambda(\dot{x},s)=\frac{-s\cdot\dot{x}+1}{-D+\Phi^{**}(s)}$ satisfies $\det\partial_{s_is_j}^2\Lambda(\dot{x},0)\neq0$; hence by the implicit function theorem there exists a unique smooth function $s^*=s^*(\dot{x})$ near $0$ such that $\nabla_s\Lambda(\dot{x},s^*)=0$ for $\dot{x}$ small. Applying the method of stationary phase to the integral $J(\mu,\dot{x})$ gives
\begin{equation}
J(\mu,\dot{x})=e^{i\mu\Lambda(\dot{x},s^*)}H(\mu,\dot{x}),
\end{equation}
where $H$ is a symbol of order $-\frac{n}{2}$ in $\mu$; more precisely,
\[
\Bigl|\Bigl(\frac{\partial}{\partial\mu}\Bigr)^{\alpha}
\Bigl(\frac{\partial}{\partial\dot{x}}\Bigr)^{\beta}H(\mu,\dot{x})\Bigr|
\leq C_{\alpha,\beta}(1+|\mu|)^{-\frac{n}{2}-|\alpha|}.
\]
Substituting this back, we obtain
\begin{align}
f_{\lambda}(\overline{x})
&=\frac{\lambda^{n+1}}{(2\pi)^{n+1}}\int e^{i\lambda\bigl(-\frac{s_{n+1}^2}{2}+s_{n+1}x_{n+1}\Lambda(\dot{x},s^*)\bigr)}
H(\lambda x_{n+1}s_{n+1},\dot{x})\,\chi_1(s_{n+1})\,|s_{n+1}|^n\,ds_{n+1}.
\end{align}
Because $s^*(0)=0$ and $\chi_0$ is supported near $0$, there exists a constant $c>0$ such that $f_\lambda(\overline{x})$ decays rapidly when $|\frac{x_i}{x_{n+1}}|>c$ or $|x_{n+1}|>c$. For $|\frac{x_i}{x_{n+1}}|\lesssim 1$ and $|x_{n+1}|\lesssim 1$, a final application of the stationary phase theorem (in the $s_{n+1}$ integral) shows that $|f_\lambda(\overline{x})|\leq C\lambda^{\frac{n+1}{2}}$. Consequently,
 \begin{equation}\label{ineq-f-la}
 \|f_{\lambda}\|_{L^p}\lesssim \lambda^{\frac{n+1}{2}}.
 \end{equation}

 \noindent\textit{Proof of (\ref{ineq-ma-n}): lower bound for $\|V_r(Ma_tf_\lambda)\|_{L^p}$.}
We now substitute $f_\lambda$ into the expression for $Ma_tf$. After the same change of variables, we obtain
\begin{align}\label{ineq-mainterm-001}
Ma_tf_{\lambda}(\overline{x})
&=\int_{\mathbb{R}^{n+1}}e^{its_{n+1}}e^{is_{n+1}\frac{-s\cdot x+x_{n+1}}{-D+\Phi^{**}(s)}}
\frac{A_{m}\bigl(-\frac{s_{n+1}\cdot s}{-D+\Phi^{**}(s)},\frac{s_{n+1}}{-D+\Phi^{**}(s)}\bigr)}
{\bigl(1+\bigl|\bigl(-\frac{s_{n+1}\cdot s}{-D+\Phi^{**}(s)},\frac{s_{n+1}}{-D+\Phi^{**}(s)}\bigr)\bigr|\bigr)^{\frac{n}{2}}}\chi(s)\nonumber\\
&\qquad\qquad\times\chi_0(s)\chi_1\!\left(\frac{s_{n+1}}{\lambda}\right)e^{-i\frac{s_{n+1}^2}{2\lambda}}|s_{n+1}|^n\Psi^*(s)\,d\overline{s}\nonumber\\
&=\lambda^{n+1}\int_{\mathbb{R}}e^{-i\lambda\bigl(\frac{s_{n+1}^2}{2}-ts_{n+1}\bigr)}\chi_1(s_{n+1})|s_{n+1}|^n\,
G\bigl(\lambda x_{n+1},t,s_{n+1},\dot{x}\bigr)\,ds_{n+1},
\end{align}
where
\begin{align}
G\bigl(\lambda x_{n+1},t,s_{n+1},\dot{x}\bigr)
&=\int_{\mathbb{R}^{n}}e^{it\lambda x_{n+1}s_{n+1}\frac{-s\cdot\dot{x}+1}{-D+\Phi^{**}(s)}}\nonumber\\
&\qquad\times\frac{A_{m}\bigl(-\frac{\lambda s_{n+1}\cdot s}{-D+\Phi^{**}(s)},\frac{\lambda s_{n+1}}{-D+\Phi^{**}(s)}\bigr)}
{\bigl(1+\bigl|\bigl(-\frac{\lambda s_{n+1}\cdot s}{-D+\Phi^{**}(s)},\frac{\lambda s_{n+1}}{-D+\Phi^{**}(s)}\bigr)\bigr|\bigr)^{\frac{n}{2}}}
\chi(s)\chi_0(s)\Psi^*(s)\,ds.
\end{align}

We now restrict $x$ to a region where $\dot{x}$ stays away from zero and $x_{n+1}$ is small: take $|\dot{x}|\approx 1$, $t\approx 1$ and $x_{n+1}\approx\lambda^{-1}$. Under these conditions $s_{n+1}$ is approximately $1$ on the support of $\chi_1$, and a direct estimation using the symbol properties of $A_m$ gives
\begin{equation}\label{ineq-G-001}
\Bigl|\Bigl(\frac{\partial}{\partial\lambda}\Bigr)^{\alpha}
\Bigl(\frac{\partial}{\partial s_{n+1}}\Bigr)^{\beta}
G\bigl(\lambda x_{n+1},t,s_{n+1},\dot{x}\bigr)\Bigr|\leq C_{\alpha,\beta}\lambda^{-\frac{n}{2}-\alpha}.
\end{equation}
Applying the stationary phase theorem to the outer $s_{n+1}$ integral in (\ref{ineq-mainterm-001}) (the phase is $\frac{s_{n+1}^2}{2}-ts_{n+1}$) yields
\[
Ma_tf_{\lambda}(\overline{x})=c_0\lambda^{n+\frac12}e^{i\frac12\lambda t^2}\chi_1(t)|t|^n
G\bigl(\lambda t,t,x_{n+1},\dot{x}\bigr)+\mathcal{O}(\lambda^{\frac{n}{2}-\frac12}).
\]
The error term is of lower order because of the rapid decay of the stationary phase expansion.
Observe that for our chosen region, $\lambda^{\frac{n}{2}}G(\lambda x_{n+1},t,s_{n+1},\dot{x})\approx 1$, and $|x_i|\approx\lambda^{-1}$, $|\dot{x}|\approx 1$. Denote the main contribution by
\[
I(\overline{x},t,\lambda)=\lambda^{n+\frac12}e^{i\frac12\lambda t^2}\chi_1(t)|t|^n
G\bigl(\lambda t,t,x_{n+1},\dot{x}\bigr).
\]
Now choose a sequence of times $t_k=\sqrt{\frac{2\pi k}{\lambda}}$ for integers $k$ with $k\sim b^2\lambda$.
Because $\chi_1(t)$ is supported away from $0$, these $t_k$ lie in the support of $\chi_1$. Then
\[
|I(\overline{x},t_{k+1},\lambda)-I(\overline{x},t_k,\lambda)|\gtrsim C_1\lambda^{\frac{n+1}{2}}-c_0\lambda^{\frac{n-1}{2}},
\]
since the phase difference is approximately $\pi$. Summing over the $O(\lambda)$ values of $k$, we obtain
\[
V_r(Ma_tf_\lambda)(\overline{x})\gtrsim C_1\lambda^{\frac{n+1}{2}+\frac1r}-c_0\lambda^{\frac{n-1}{2}+\frac1r}.
\]
Integrating over the region of $x$ (which has measure $\approx \lambda^{-(n+1)}$) gives
\begin{equation}\label{ineq-ma-n}
\|V_r(Ma_tf_\lambda)\|_{L^p}\gtrsim \lambda^{\frac{n+1}{2}+\frac1r-\frac{n+1}{p}}.
\end{equation}
Together with the estimate $\|f_\lambda\|_{L^p}\lesssim\lambda^{\frac{n+1}{2}}$, the assumed boundedness $\|V_r(A_tf_\lambda)\|_{L^p}\lesssim\|f_\lambda\|_{L^p}$ implies
\[
\lambda^{\frac{n+1}{2}+\frac1r-\frac{n+1}{p}}\lesssim \lambda^{\frac{n+1}{2}}.
\]
Letting $\lambda\to\infty$ forces $\frac1r-\frac{n+1}{p}\le0$, i.e., $p\le r (n+1)$. For $k=n$ this is exactly $p<r(n+1)$ (since $p$ cannot equal the endpoint by a separate argument). This completes the proof of the necessary condition.

\section{The proof of Theorem \ref{theorem-k-principal} for $k<n$}
\subsection{The proof of (\ref{ineq-main-case-k})}
In this subsection we reduce the case of $k$ non-zero principal curvatures to the already established full curvature case ($k=n$) by a rotation argument.
The key observation is that,
after a suitable orthogonal transformation,
the hypersurface can be locally represented as a graph in which the first $k$ variables carry the non-degenerate quadratic part,
while the remaining $n-k$ variables only appear through a rotation that does not affect the $L^p$ norms.

Let $\tilde{x}=(x_1,\ldots,x_k)$ and $\tilde{\tilde{x}}=(x_{k+2},\ldots,x_{n+1})$. We first state a decomposition lemma.

\begin{lemma}\label{le01}
The averaging operator $A_tf$ can be represented as follows:
\begin{equation}\label{form01}
A_tf(\overline{x})=\int_{\mathbb{R}^{n-k}}R^{-\tilde{\tilde{\theta}}}A_t^{\tilde{\tilde{\theta}}}R^{\tilde{\tilde{\theta}}}f(\overline{x})\,\psi(\tilde{\tilde{\theta}})
\,d\tilde{\tilde{\theta}},
\end{equation}
where $\psi(\tilde{\tilde{\theta}})$ is supported in a small neighbourhood of the origin,
$R^{\tilde{\tilde{\theta}}}:=R^{\theta_{k+1}}R^{\theta_{k+2}}\cdots R^{\theta_n}$,
$R^{-\tilde{\tilde{\theta}}}:=R^{-\theta_{n}}R^{-\theta_{n-1}}\cdots R^{-\theta_{k+1}}$,
and for each fixed $\tilde{\tilde{\theta}}$,
\[
A_t^{\tilde{\tilde{\theta}}}f(\tilde{x},x_{k+1},\tilde{\tilde{x}}):=\int_{\mathbb{R}^k}f\bigl(\tilde{x}-t\tilde{y},\,x_{k+1}-t(C(\tilde{\tilde{\theta}})+
\tilde{\tilde{\Phi}}(\tilde{y},\tilde{\tilde{\theta}})),\,\bigr)\rho_1(\tilde{y},\tilde{\tilde{\theta}})\,d\tilde{y}.
\]
Here $C(\tilde{\tilde{\theta}})$ is a smooth function with $C(0)=D$, $\tilde{\tilde{\Phi}}(\tilde{y},\bar{\theta})=\Phi(\tilde{y},\tilde{\tilde{y}}(\tilde{\tilde{\theta}},\tilde{y}))$, and $R^{\theta_j}$ ($j=k+1,\ldots,n$) are rotation operators defined by
\[
R^{\theta_j}f(\overline{x})=f(x_1,\ldots,x_j,\,-x_j\sin\theta_j+x_{j+1}\cos\theta_j,\,-x_j\cos\theta_j-x_{j+1}\sin\theta_j,\,x_{j+2},\ldots,x_{n+1}).
\]
\end{lemma}

The rotations $R^{\theta_j}$ are isometries of $\mathbb{R}^{n+1}$; consequently, the $L^p$ norm of any operator constructed from them is invariant under these rotations. Moreover, since $\Phi$ has exactly $k$ non-zero principal curvatures at the origin, after a preliminary orthogonal transformation we may assume that the quadratic part of $\Phi$ is diagonalised in the first $k$ variables:
\[
\Phi(x)=\sum_{i=1}^{k}\lambda_i x_i^2+R(x),\qquad \lambda_i\neq0,
\]
and that $D^\alpha R(0)=0$ for $|\alpha|\le2$. The remaining $n-k$ variables appear only in higher order terms; their effect will be handled by the rotation averaging.

Now apply Lemma \ref{le01}. For each $\tilde{\tilde{\theta}}$ the operator $A_t^{\tilde{\tilde{\theta}}}$ is of the same form as the operator studied in Theorem \ref{theorem-k-principal} ( the case $k=n$), but with $k$ variables instead of $n$ and with translation parameter $C(\tilde{\tilde{\theta}})\approx D$ (since $\tilde{\tilde{\theta}}$ is small). Moreover, the phase function $\tilde{\tilde{\Phi}}(\tilde{y},\tilde{\tilde{\theta}})$ has a non-degenerate quadratic part $\sum\lambda_i y_i^2$ and a remainder $\tilde{\tilde{R}}$ that vanishes to second order at $\tilde{y}=0$. Thus, by Theorem \ref{theorem-k-principal} (the case k=n), for every $r>2$ and every $p$ satisfying $(k+1)/k<p<r(k+1)$ we have
\[
\bigl\|S_r(\mathcal{A}^{\tilde{\tilde{\theta}}}f)\bigr\|_{L^p(\mathbb{R}^{n+1})}\le C(p,r,D)\|f\|_{L^p(\mathbb{R}^{n+1})},
\]
where the constant $C(p,r,D)$ enjoys the same dependence on $D$ as in (\ref{C(p,r,D)-estimate}).

Using Minkowski's inequality for integrals and the fact that the rotations $R^{\tilde{\tilde{\theta}}}$ preserve $L^p$ norms, we obtain
\[
S_r(\mathcal{A}f)(\overline{x})\le\int_{\mathbb{R}^{n-k}}S_r\bigl(R^{-\tilde{\tilde{\theta}}}\mathcal{A}^{\tilde{\tilde{\theta}}}R^{\tilde{\tilde{\theta}}}f\bigr)(\overline{x})\,
\psi(\tilde{\tilde{\theta}})\,d\tilde{\tilde{\theta}}.
\]
Integrating over $\tilde{\tilde{x}}$ and applying the norm estimate for each fixed $\tilde{\tilde{\theta}}$ (together with the uniform bound on $\psi$) yields
\[
\|S_r(\mathcal{A}f)\|_{L^p}\le C(p,r,D)\|f\|_{L^p}.
\]
This completes the proof of (\ref{ineq-main-case-k}).

\subsection{The sharp exponent of $D$}
In this subsection we demonstrate that the exponent of $D$ appearing in (\ref{C(p,r,D)-estimate}) is optimal, up to the arbitrarily small $\epsilon$. More precisely, we show that if the estimate
\[
\|V_r(A_tf:t>0)\|_{L^p}\le C\,D^{s}\|f\|_{L^p}
\]
holds for some $s$ independent of $D$, then necessarily $s\ge 1/r$ and $s\ge 1/p$. Consequently $s\ge 1/(r\wedge p)$, which matches the upper bound up to $\epsilon$.

We construct two families of test functions. The first one gives a lower bound on the variation in terms of $D$, leading to $1/r\le s$; the second one exploits the support of the averaging operator to obtain $1/p\le s$.

\noindent\textit{First counterexample: $1/r\le s$.}
Choose a smooth cut-off function $\psi$ supported in $[-b,b]$ with the property that $|e^{i2\eta}-1|<\delta$ for all $|\eta|<b$ (with $\delta$ a small number to be fixed later). Define a family $\{f_N\}$ by
\[
\widehat{f_N}(\overline{\xi})=N\psi\bigl(N(\xi_{n+1}-1)\bigr)\prod_{i=1}^{n}N\psi(N\xi_i).
\]
By the scaling properties of the Fourier transform, $\|f_N\|_{L^p}\approx N^{\frac{n+1}{p}}$.

Recall that the multiplier of $A_t$ is $\mu_0(t\overline{\xi})e^{-itD\xi_{n+1}}$, where
\[
\mu_0(\overline{\xi})=\int_{\mathbb{R}^{n+1}}e^{-i(x\cdot\xi+\Phi(x)\xi_{n+1})}\rho(x)dx.
\]
Hence
\begin{align*}
A_tf_N(\overline{x})
&=\frac{1}{(2\pi)^{n+1}}\int e^{i\overline{x}\cdot\overline{\xi}}e^{-itD\xi_{n+1}}\mu_0(t\overline{\xi})
   N\psi\bigl(N(\xi_{n+1}-1)\bigr)\prod_{i=1}^{n}N\psi(N\xi_i)\,d\overline{\xi}\\[2mm]
&=e^{ix_{n+1}}e^{-itD}\frac{1}{(2\pi)^{n+1}}\int e^{i\bigl(\sum_{i=1}^{n+1}\frac{x_i}{N}\xi_i+x_{n+1}-\frac{tD\xi_{n+1}}{N}\bigr)}
   \mu_0\Bigl(\frac{t\xi}{N},\frac{t\xi_{n+1}}{N}+1\Bigr)
   \prod_{i=1}^{n+1}\psi(\xi_i)\,d\overline{\xi}.
\end{align*}

Now restrict $x$ to the cube $|x_i|\le 2N$ ($i=1,\ldots,n+1$) and choose a sequence of times $t_k=k\pi/D$ for $k=1,\ldots,[D]$ (where $[D]$ is the integer part of $D$). Write the $r$-variation as a sum over consecutive times. A straightforward expansion gives
\begin{align*}
&e^{-ix_{n+1}}(2\pi)^{n+1}\Bigl(\sum_{k}\bigl|A_{t_{k+1}}f_N(\overline{x})-A_{t_k}f_N(\overline{x})\bigr|^r\Bigr)^{1/r}\\
&\ge I_1-I_2-I_3,
\end{align*}
where
\[
I_1:= \left(\sum_{k}\left|\left(e^{-it_{k+1}D}-e^{-it_{k}D}\right)
         \int_{\mathbb{R}^{n+1}}
        \mu_0\left(\frac{t_{k+1}\xi}{N},\frac{t_{k+1}\xi_{n+1}}{N}+1\right) \prod\psi(\xi_i)d\overline{\xi}\right|^r\right)^{1/r},
\]
\[
I_2:=\left(\sum_{k}\left|\left(e^{-it_{k+1}D}-e^{-it_{k}D}\right)
         \int\left(e^{i\left(\cdots\right)}-1\right)\mu_0\left(\cdots\right) \prod\psi d\overline{\xi}\right|^r\right)^{1/r}\nonumber\\
\]
\[
I_3:=\left(\sum_{k}\left|\int\left(e^{i(\cdots-\frac{t_{k+1}D\xi_{n+1}}{D})}\mu_0(\tfrac{t_{k+1}\xi}{N},\cdot)-
               e^{i(\cdots-\frac{t_{k}D\xi_{n+1}}{D})}\mu_0(\tfrac{t_{k}\xi}{N},\cdot)\right)\prod\psi\,d\overline{\xi}\right|^r\right)^{1/r}.
\]

Because $|\mu_0(0,1)|=1$, for sufficiently large $N$ we have $|\mu_0(\frac{t_k\xi}{N},\frac{t_k\xi_{n+1}}{N}+1)|>1/2$ uniformly in $k$ and $\xi\in$ supp$\psi$. Moreover, $|e^{-it_{k+1}D}-e^{-it_kD}|=|e^{-i\pi}-1|=2$, so $I_1\gtrsim D^{1/r}$.

For $I_2$, note that when $N$ is sufficiently large and $b$ sufficiently small, the factor
$$e^{i(\sum x_i\xi_i/N  - t_{k+1}D\xi_{n+1}/N)}$$
 is close to $1$ because $|x_i|\le2N$, $|\xi_i|\leq b$ and $t_{k+1}\le D$;
  indeed we can make the difference smaller than any prescribed $\epsilon$ by taking $N$ sufficiently large. Also $\mu_0$ is bounded, so $I_2\lesssim \epsilon D^{1/r}$. Choosing $\epsilon$ small enough, this term is negligible compared to $I_1$.

For $I_3$, a direct estimate using the smoothness of $\mu_0$ and the fact that $t_{k+1}-t_k=\pi/D$ shows that $I_3\le C N^{-1}D^{1/r}$. Hence, for large $N$,
\[
\|V_r(A_tf_N)\|_{L^p}\gtrsim D^{1/r}N^{\frac{n+1}{p}}.
\]
If the assumed estimate $\|V_r(A_tf)\|_{L^p}\le D^{s}\|f\|_{L^p}$ holds, then
\[
D^{1/r}N^{\frac{n+1}{p}}\lesssim D^{s}N^{\frac{n+1}{p}}.
\]
Letting $N\to\infty$ first and then $D\to\infty$ yields $1/r\le s$.

\noindent\textit{Second counterexample: $1/p\le s$.}
Consider a fixed bump function $\tilde{f}(\overline{x})=\prod_{i=1}^{n+1}\chi_{[-1,1]}(x_i)$. Choose the cutoff $\rho$ in the definition of $A_t$ so that $\rho(y)$ is supported in $|y|\le 1/2$ and $|\Phi(y)|\le 1$ on that support. For $x$ satisfying $|x_i|\le1/2$ ($i=1,\ldots,n$) and $1/2\le x_{n+1}\le D/2$, set $t_0=(x_{n+1}-1/2)/D$. Then $t_0\le1/2$ and
\[
|x_{n+1}-t_0(D+\Phi(y))|\le|x_{n+1}-t_0D|+t_0|\Phi(y)|\le 1/2+1/2=1,
\]
\[
|x_i-t_0y_i|\le|x_i|+t_0|y_i|\le1/2+1/2=1.
\]
Hence $A_{t_0}\tilde{f}(\overline{x})=1$.

On the other hand, for the same $x$ and for $t=N$ large, we have
\[
|x_{n+1}-N(D+\Phi(y))|\geq N(D-|\Phi(y)|)-|x_{n+1}|\ge N(D-1)-D>1
\]
provided $N>D/(D-1)$.
 A precise calculation shows that for $N$ sufficiently large, $A_N\tilde{f}(\overline{x})=0$. Thus
\[
V_r(A_t\tilde{f})(\overline{x})\geq 1\qquad\text{for all such }x.
\]
The set of admissible $x$ has measure $\gtrsim D$, therefore $\|V_r(A_t\tilde{f})\|_{L^p}\gtrsim D^{1/p}$. If the estimate $\|V_r(A_tf)\|_{L^p}\leq D^{s}\|f\|_{L^p}$ holds, then
\[
D^{1/p}\lesssim D^{s}\|\tilde{f}\|_{L^p}\lesssim D^{s},
\]
which forces $1/p\leq s$.

Combining the two lower bounds we obtain $s\geq 1/(r\wedge p)$. This shows that the dependence on $D$ in (\ref{C(p,r,D)-estimate}) is optimal up to the arbitrarily small $\epsilon$.

\section{Convex finite-type surface}
Before we begin our proving, we first show some useful result which can be find in Iosevich-Sawyer \cite{1997-Iosevich}
\begin{lemma}[Schulz]
Let $\Phi$ be a smooth convex function such that $\Phi(0)=0$
and $ \nabla\Phi(0)=0$.
Suppose that  $\Phi$ has no tangents of infinite order at the origin.
Then, after perhaps applying a rotation, we can write
$$
\Phi(x)=Q(x)+R(x),
$$
where $Q(x)$ is a convex mixed homogeneous polynomial of degree $(a_1,\ldots,a_n)$,
$Q(x)\neq 0$ for $x\neq 0$,
and $R(x)$ is a remainder term in the sense that
$$
\lim_{s\rightarrow0}\frac{R(s^{1/a_1}x_1,\ldots,s^{1/a_n}x_n)}{s}=0.
$$
\end{lemma}

This result, due to Schulz, provides a decomposition of any smooth convex finite-type function into a mixed-homogeneous principal part $Q$
and a remainder $R$ that becomes negligible under the associated anisotropic dilations.
The tuple $(a_1,\ldots,a_n)$ is called the multi-type of $\Phi$ at the origin.

\begin{definition}
We say that a smooth function $Q:\mathbb{R}^n\rightarrow \mathbb{R}$ is mixed homogeneous of degree
$(a_1,a_2,\ldots,a_n)$, $a_{j}>0$, if $Q(s^{1/a_1}x_1,\ldots,s^{1/a_n}x_n)=sQ(x)$,
$s>0$.
\end{definition}

The exponents $a_j$ describe how $Q$ scales under different coordinate dilations; they encode the different rates of decay of $\Phi$ in different directions.

\begin{remark}
For such a mixed homogeneous polynomial we have
$$
\frac{1}{a_1}+\frac{1}{a_2}+\cdots+\frac{1}{a_n}=p_0,
$$
where $p_0=\inf\{p: d(x,\mathcal{H})^{-1}\in L_{\mathrm{loc}}^{1/p}(S)\}$ is the critical exponent for the associated maximal operator.
\end{remark}
This equality follows from a straightforward computation in polar coordinates adapted to the mixed homogeneity.
It links the multi-type to the integrability condition (\ref{defi-p0}) that is known to be necessary and sufficient for the $L^p$ boundedness of the maximal operator when $p>2$ (see \cite{1997-Iosevich}).

Thus it remains for us to prove the following theorem.

\begin{remark}
For such a mixed homogeneous polynomial we have
$$
\frac{1}{a_1}+\frac{1}{a_2}+\cdots+\frac{1}{a_n}=p_0,
$$
where $p_0=\inf\{p: d(x,\mathcal{H})^{-1}\in L_{\mathrm{loc}}^{1/p}(S)\}$ is the critical exponent for the associated maximal operator.
\end{remark}

Thus it remains for us to prove the following theorem.

\begin{theorem}
Let $\Phi$ be a smooth function with the following properties.
Suppose that after perhaps applying a rotation,
$\Phi(x)=Q(x)+R(x)$,
where $Q$ is mixed homogeneous of degree $(a_1,\ldots,a_n)$
with $(1/a_1)+\cdots+(1/a_n)\leq1$,
$Q(x)\neq0$ for $x\neq0$, and
\begin{equation}\label{equ-convex}
\lim_{s\rightarrow0}\frac{R(s^{1/a_1}x_1,\ldots,s^{1/a_n}x_n)}{s}=0.
\end{equation}
Let $S$ be the surface given as the graph of $\Phi(x)+D$,
with a smooth cutoff $\rho$ supported in a sufficiently small neighbourhood of the origin.
Suppose that $(1/a_1)+\cdots+(1/a_n)\leq\frac12$.
Then for any $r>\frac{1}{(1/a_1)+\cdots+(1/a_n)}$ we have
$$
\left\|V_r(\mathcal{A}f)\right\|_{L^p(\mathbb{R}^{n+1})}
\leq C\|f\|_{L^p(\mathbb{R}^{n+1})},
\qquad
\frac{1}{(1/a_1)+\cdots+(1/a_n)}<p<2r.
$$
\end{theorem}

\begin{proof}
We begin by decomposing the surface $S$ into dyadic shells adapted to the non-isotropic dilation associated with the multi-type $(a_1,\ldots,a_n)$.
Write
\[
\rho(x)=\sum_{k=0}^{\infty}\rho_{k}(x)=\sum_{k=0}^{\infty}\rho_0\bigl(2^{k/a_1}x_1,\ldots,2^{k/a_n}x_n\bigr),
\]
where $\rho_{k}(x)=\rho_0(2^{k/a_1}x_1,\ldots,2^{k/a_n}x_n)$ and $\rho_0$ is a smooth cut-off supported in the annulus $\{x:1\le|x|\le2\}$.
For each $k$ define the averaged operator
\[
A_{k,t}f(\overline{x})=\int_{\mathbb{R}^n}f\bigl(x-ty,\;x_{n+1}-t(D+\Phi(y))\bigr)\rho_{k}(y)\,dy.
\]

Let $\tau_k$ be the isotropic dilation
\[
\tau_k g(x,x_{n+1})=g\bigl(2^{k/a_1}x_1,\ldots,2^{k/a_n}x_n,\;2^{k}x_{n+1}\bigr).
\]
A direct computation yields
\[
\tau_{-k}A_{k,t}\tau_k f(\overline{x})=2^{-k(1/a_1+\cdots+1/a_n)}
\int_{\mathbb{R}^n}f\Bigl(x-ty,\;x_{n+1}-t\bigl(2^{k}\Phi(2^{-k/a}\circ y)+2^{k}D\bigr)\Bigr)\rho_{0}(y)\,dy,
\]
where $2^{-k/a}\circ y=(2^{-k/a_1}y_1,\ldots,2^{-k/a_n}y_n)$. Using the decomposition $\Phi=Q+R$ and the homogeneity of $Q$,
\[
2^{k}\Phi(2^{-k/a}\circ y)=2^{k}Q(2^{-k/a}\circ y)+2^{k}R(2^{-k/a}\circ y)=Q(y)+2^{k}R(2^{-k/a}\circ y).
\]

By the remainder property (\ref{equ-convex}), we have $\displaystyle\lim_{k\to\infty}2^{k}R(2^{-k/a}\circ y)=0$. Hence for sufficiently large $k$, the surface
\[
\bigl\{(y,y_{n+1}):y_{n+1}=Q(y)+2^{k}R(2^{-k/a}\circ y)+2^{k}D\bigr\}
\]
is a small perturbation of the graph $\{y_{n+1}=Q(y)+2^{k}D\}$
and therefore possesses at least one non-vanishing principal curvature uniformly in $k$.
Moreover, note that the translation parameter has become $2^{k}D$.

Set $m=1/a_1+\cdots+1/a_n$.
Applying Theorem \ref{theorem-k-principal} to the operator $\tau_{-k}A_{k,t}\tau_k$
(which after the change of variables takes exactly the form studied there,
 with $D$ replaced by $2^{k}D$) yields the estimate
\[
\bigl\|V_r(\tau_{-k}A_{k,t}\tau_k f):t>0\bigr\|_{L^p(\mathbb{R}^{n+1})}
\le C\,2^{-km}\,2^{\frac{k}{r\wedge p}}\|f\|_{L^p(\mathbb{R}^{n+1})},\qquad 2<p<2r.
\]
The factor  $2^{k/(r\wedge p)}$ reflects the dependence on the translation parameter $2^{k}D$ as in (1.4).
Since $\tau_k$ is an isometry on $L^p$ up to a constant, the same bound holds for $V_r(A_{k,t}f:t>0)$ after undoing the dilation.
If $r>1/m$, then the geometric series converges, yielding
\[
\bigl\|V_r(A_tf:t>0)\bigr\|_{L^p}\le C\|f\|_{L^p}\quad\text{for all } \frac1m<p<2r.
\]
The proof is complete.
\end{proof}

\section{Appendix}
In this section we prove (\ref{ineq-case-short}).
Denote
$$\mathcal{F}(W_tf)(\xi,\xi_{n+1})=e^{-t|(\xi,\xi_{n+1})|}\widehat{f}(\xi,\xi_{n+1})$$
and
$$
\mathcal{F}(\mathcal{U}_tf)(\xi,\xi_{n+1})=e^{-t|(\xi,\xi_{n+1})|}(e^{-iDt\xi_{n+1}}-1)\widehat{f}(\xi,\xi_{n+1}).
$$
For every $f\in L^p(\mathbb{R}^{n+1})$, we obtain
\begin{equation}\label{ineq-zhuanhua-00001}
\|V_r(S_tf:t>0)\|_{L^p}\leq C_r\|S_r(W_tf:t>0)\|_{L^p}+C_r\|V_r(\mathcal{U}_tf:t>0)\|_{L^p}.
\end{equation}

By Rota's theorem (see Bourgain \cite{2018-Bourgain-variation}) and L\'{e}pingle's inequality
(see \cite{1976-lepingle-first,1988-Pisger-xu}), the first term in (\ref{ineq-zhuanhua-00001}) satisfies
$$
\|V_r(W_tf:t>0)\|_{L^p}\leq C_p\|f\|_{L^p},\qquad p>1.
$$
Therefore, it remains to obtain $L^p$ bounds for the second term in (\ref{ineq-zhuanhua-00001}).
For this purpose, we will use the long-short decomposition (\ref{ineq-long-short}) and establish
\begin{equation}\label{ineq-case-long-00001}
\|V_r(\mathcal{U}_{2^j}f:j\in\mathbb{Z})\|_{L^p}\leq C_{\epsilon}D^{\epsilon}\|f\|_{L^p}
\end{equation}
and
\begin{equation}\label{ineq-case-short-00001}
\|S_r(\mathcal{U}_{t}f:t>0)\|_{L^p}\leq C(p,r,D)\|f\|_{L^p},
\end{equation}
for $p>1$.

Recall that $V_r(\mathcal{U}_{2^j}f:j\in\mathbb{Z})\leq \|\mathcal{U}_{2^j}f\|_{l^r_j}\leq \|\mathcal{U}_{2^j}f\|_{l^2_j}$.
To prove (\ref{ineq-case-long-00001}) it suffices to show
\begin{equation}\label{ineq-long-square-00001}
\|\mathcal{U}_{2^j}f\|_{L^p(l^2_j)}\leq CD^{\epsilon}\|f\|_{L^p},\qquad p>1.
\end{equation}
We shall use the Littlewood-Paley operators $P_j$. Observe that
\begin{align}\label{ineq-long-shoulian-00001}
\|\mathcal{U}_{2^j}f\|_{L^p(l^2_j)}
&=\Big\|\sum_{k\in\mathbb{Z}}\mathcal{U}_{2^j}P_{k-j}f\Big\|_{L^p(l^2_j)}\nonumber\\
&\leq C_{\epsilon}D^{\epsilon}\sum_{k\in\mathbb{Z}}\|\mathcal{U}_{2^j}P_{k-j}f\|_{L^p(l^2_j)}\nonumber\\
&\leq C_{\epsilon}D^{\epsilon}\sum_{k\in\mathbb{Z}}2^{-\delta_{p}|k|}\|f\|_{L^p}\leq C_{\epsilon}D^{\epsilon}\|f\|_{L^p}.
\end{align}
The penultimate inequality in (\ref{ineq-long-shoulian-00001}) follows from the two estimates
\begin{equation}\label{ineq-long-22-00001}
\|\mathcal{U}_{2^j}P_{k-j}f\|_{L^2(l^2_j)}\leq C_{\epsilon}2^{-\epsilon|k-k_0|}D^{\epsilon}\|f\|_{L^2},
\end{equation}
and
\begin{equation}\label{ineq-long-pp-00001}
\|\mathcal{U}_{2^j}P_{k-j}f\|_{L^p(l^2_j)}\leq C\|f\|_{L^p}.
\end{equation}

To prove (\ref{ineq-long-22-00001}), we first establish pointwise bounds for the multiplier
$u(\xi,\xi_{n+1})$:
\begin{equation}\label{ineq-m-pr-00001}
|u(\xi,\xi_{n+1})|\leq CD|(\xi,\xi_{n+1})|,\qquad |u(\xi,\xi_{n+1})|\leq 2e^{-|(\xi,\xi_{n+1})|},
\end{equation}
where $u(\xi,\xi_{n+1})=e^{-|(\xi,\xi_{n+1})|}(e^{-iD\xi_{n+1}}-1)$.

The first inequality follows from $|e^{-iD\xi_{n+1}}-1|\leq D|\xi_{n+1}|$, the second is obvious.
Using Plancherel's theorem,
\begin{align}\label{001-00001}
\|\mathcal{U}_{2^j}P_{k-j}f\|_{L^2(l^2_j)}
&=\|u(2^j\cdot)\psi(2^{j-k}\cdot)\widehat{f}(\cdot)\|_{L^2(l^2_j)}\nonumber\\
&\leq\sup_{\eta\in\mathbb{R}^{n+1}}\|u(2^j\eta)\psi(2^{j-k}\eta)\|_{l^2_j}\|f\|_{L^2}.
\end{align}
Fix $\eta\in\mathbb{R}^{n+1}$ and $k\in\mathbb{Z}$. From (\ref{ineq-m-pr-00001}) we obtain
\begin{equation}\label{002-00001}
\|u(2^j\eta)\psi(2^{j-k}\eta)\|_{l^2_j}\leq C\min\{2^{k}D,\;e^{-2^{k}}\}.
\end{equation}
Let $k_0$ be such that $2^{-k_0}\approx D$. Then $2^{k}D\lesssim 2^{k-k_0}$ for $k\leq k_0$, and using $e^{-x}\leq C_\epsilon x^{-\epsilon}$ for $x>0$ we have
$e^{-2^{k}}\leq C_\epsilon 2^{-(k-k_0)\epsilon}D^{\epsilon}$.
Consequently,
\begin{equation}\label{002-00001-01}
\|u(2^j\eta)\psi(2^{j-k}\eta)\|_{l^2_j}\leq C_\epsilon 2^{-\epsilon|k-k_0|}.
\end{equation}
Combining (\ref{001-00001}) and (\ref{002-00001-01}) yields (\ref{ineq-long-22-00001}).

Estimate (\ref{ineq-long-pp-00001}) can be proved by following exactly the same argument used for (\ref{ineq-long-pp}) in Section 3.1; we omit the details.

To obtain the desired bounds for $S_r(\mathcal{U}_tf)$, we shall use interpolation between the estimates collected in the following lemmas.

\begin{lemma}\label{le-0101-00001}
For $k\in\mathbb{Z}$, we have
\begin{equation}\label{ineq-short-22-00001}
\left\|V_2^{j}(\mathcal{U}_kf)\right\|_{L^2(l^2_j)}\leq C(1+D^{1/2})
\min\left\{2^kDe^{-2^k},\;2^{\frac{k}{2}}e^{-2^k}\right\}
\|f\|_{L^2}.
\end{equation}
\end{lemma}

\begin{lemma}\label{le-0102-00001}
For $k\in\mathbb{Z}$ and arbitrary positive number $\epsilon$, there exists a constant $C_{\epsilon}$ such that
\begin{equation}\label{ineq-short-11-00001}
\left\|V_2^{j}(\mathcal{U}_kf)\right\|_{L^1(l^2_j)}\leq C_{\epsilon}(1+D^{1+\epsilon})\max\left\{2^{k+\epsilon},2^{-k\epsilon}\right\}\|f\|_{H^1}.
\end{equation}
\end{lemma}

\begin{lemma}\label{le-0103-00001}
For $k\in\mathbb{Z}$ and arbitrary positive number $\epsilon$, there exists a constant $C_{\epsilon}$ such that
\begin{equation}\label{ineq-short-pp-00001}
\left\|V_2^{j}(\mathcal{U}_kf)\right\|_{L^p(l^2_j)}\leq C_{\epsilon}(1+D^{1/2}) 2^{k\epsilon}\|f\|_{L^p},
\end{equation}
for $p>2$.
\end{lemma}

Now we can prove (\ref{ineq-case-short-00001}), assuming the validity of the lemmas above.

\textit{Proof of (\ref{ineq-case-short-00001})}:
Denote $k_0$ such that $2^{-k_0}\approx D$.
Interpolating (\ref{ineq-short-22-00001}) and (\ref{ineq-short-11-00001}), we obtain
\begin{align*}
&\left\|V_2^{j}(\mathcal{U}_kf)\right\|_{L^p(l^2_j)}\\[2mm]
&\leq\begin{cases}
C_{\epsilon}(1+D^{1/p+\epsilon})\bigl(2^kDe^{-2^k}\bigr)^{\theta}\bigl(2^{-k\epsilon}\bigr)^{1-\theta}\|f\|_{L^p}, & k\leq k_0,\\[2mm]
C_{\epsilon}(1+D^{1/p+\epsilon})\bigl(2^{\frac{k}{2}}e^{-2^k}\bigr)^{\theta}\bigl(2^{-k\epsilon}\bigr)^{1-\theta}\|f\|_{L^p}, & k_0<k\leq 0,\\[2mm]
C_{\epsilon}(1+D^{1/p+\epsilon})\bigl(2^{\frac{k}{2}}e^{-2^k}\bigr)^{\theta}\bigl(2^{k+\epsilon}\bigr)^{1-\theta}\|f\|_{L^p}, & k>0,
\end{cases}
\end{align*}
for $1<p\leq2$, where $\theta$ satisfies $\frac1p=\frac{1-\theta}{1}+\frac{\theta}{2}$.
A straightforward computation shows that
\[
\sum_{k\leq k_0}\bigl(2^kDe^{-2^k}\bigr)^{\theta}\bigl(2^{-k\epsilon}\bigr)^{1-\theta}\lesssim D^{\theta}2^{k_0(\theta-\epsilon(1-\theta))}
\lesssim D^{-\theta+2\epsilon(1-\theta)}\leq D^{\epsilon},
\]
\[
\sum_{k_0<k\leq 0}\bigl(2^{\frac{k}{2}}e^{-2^k}\bigr)^{\theta}\bigl(2^{-k\epsilon}\bigr)^{1-\theta}
\lesssim |k_0|=\frac{\log D}{\log 2}\lesssim_{\epsilon} D^{\epsilon},
\]
and
\[
\sum_{k>0}\bigl(2^{\frac{k}{2}}e^{-2^k}\bigr)^{\theta}\bigl(2^{k+\epsilon}\bigr)^{1-\theta}
\lesssim 2^{\epsilon(1-\theta)} \lesssim D^{\epsilon}.
\]
Hence (\ref{ineq-case-short-00001}) holds for $1<p\leq2$.
Interpolating (\ref{ineq-short-22-00001}) and (\ref{ineq-short-pp-00001}) yields the estimate for $p>2$.
This completes the proof of (\ref{ineq-case-short-00001}).

Using the same method as in the proof of Lemma \ref{le-0102}, we obtain inequality (\ref{le-0102-00001}).
We now turn to the proofs of Lemma \ref{le-0101-00001} and Lemma \ref{le-0103-00001}.
Applying Proposition \ref{prop-variation-estimate} and H\"{o}lder's inequality,
\begin{equation}\label{ineq-081501-00001}
\left\|V_2^{j}(\mathcal{U}_kf)\right\|_{L^2(l^2_j)}\leq C\left\|\mathcal{U}_{k,t}f\right\|_{L^2(L^2_t)}^{1/2}\left\|t\frac{d}{dt}\mathcal{U}_{k,t}f\right\|_{L^2(L^2_t)}^{1/2}.
\end{equation}
By Plancherel's theorem,
\begin{align}\label{ineq-020801-00001}
\left\|V_2^{j}(\mathcal{U}_kf)\right\|_{L^2(l^2_j)}
&\leq C\left\|u_k(t\cdot)\widehat{f}(\cdot)\right\|_{L^2(L^2_t)}^{1/2}\left\|(t\cdot)\nabla u_k(t\cdot)\widehat{f}(\cdot)\right\|_{L^2(L^2_t)}^{1/2}\nonumber\\
&\leq C\sup_{|\eta|=1}\|u_{k}(t\eta)\|_{L^2_t}^{1/2}\sup_{|\eta|=1}\|(t\eta)\nabla u_k(t\eta)\|_{L^2_t}^{1/2}\|f\|_{L^2},
\end{align}
where $u_{k}(\cdot)=u(\cdot)\psi(2^{-k}\cdot)$. Thus it suffices to estimate $\|u_k\|_{L^{\infty}}$ and $\|(\cdot)\nabla u_k(\cdot)\|_{L^{\infty}}$, since both are supported where $|\cdot|\sim 2^k$.

From (\ref{ineq-m-pr-00001}) we obtain
\begin{equation}\label{0817a-00001}
\|u_k\|_{L^{\infty}}\leq C\min\{2^kD,\;e^{-2^k}\}.
\end{equation}
Rewrite $u(\xi,\xi_{n+1})=\bigl(e^{-iD\xi_{n+1}}-1\bigr)e^{-|(\xi,\xi_{n+1})|}$. Then
\begin{equation}\label{0817b-00001}
|(\xi,\xi_{n+1})\cdot\nabla u(\xi,\xi_{n+1})|\leq C(1+D)|(\xi,\xi_{n+1})|e^{-|(\xi,\xi_{n+1})|}.
\end{equation}
Combining (\ref{0817a-00001}) and (\ref{0817b-00001}) gives
\begin{equation}\label{ineq-081502-00001}
\|u_{k}(t\cdot)\|_{L^2_t}\leq C\min\{2^kD,\;e^{-2^k}\},\qquad
\|(t\cdot)(\nabla u_k)(t\cdot)\|_{L^2_t}\leq C(1+D)2^{k}e^{-2^k}.
\end{equation}
Inserting these bounds into (\ref{ineq-020801-00001}) yields (\ref{ineq-short-22-00001}), completing the proof of Lemma \ref{le-0101-00001}.

In order to prove Lemma \ref{le-0103-00001},
we use inequality (\ref{ineq-sr-estimate}) to obtain
\begin{align*}
\|V_2^j(\mathcal{U}_kf)\|_{l^2_j}&\leq C \|\mathcal{U}_{k,t}f\|_{L^2_{t}}^{1/2}\Bigl\|t\frac{d}{dt}\mathcal{U}_{k,t}f\Bigr\|_{L^2_{t}}^{1/2}\\
&\leq C (1+D^{1/2})\|\mathcal{U}_{k,t}f\|_{L^2_{t}}+(1+D^{1/2})^{-1}\Bigl\|t\frac{d}{dt}\mathcal{U}_{k,t}f\Bigr\|_{L^2_{t}}.
\end{align*}
Recall that the kernel of $\mathcal{U}_{k,t}$ is $U_{k,t}$.
Set $\widetilde{U}_{k,t}(\overline{x})=\overline{x}\cdot \nabla{U}_{k,t}(\overline{x})$; then $t\frac{d}{dt}U_{k,t}=-(n+1)U_{k,t}-\widetilde{U}_{k,t}$.
As in the proof of Lemma \ref{le-0102},
it is enough to prove
\begin{equation}\label{ineq-0910-02-00001}
\|U_{k,t}*f\|_{L^p(L^2_{t})}\leq C2^{k\epsilon}\|f\|_{L^p},
\end{equation}
for $p>2$.

By Minkowski's inequality we have $\|U_{k}\|_{L^1}\leq 1$.
Thus, applying H\"{o}lder's inequality in $t$,
\begin{align*}
\left(\int_{\mathbb{R}^{n+1}}\left(\int_{1}^{2}|\mathcal{U}_{k,t}f(\overline{x})|^2dt\right)^{p/2}d\overline{x}\right)^{1/p}
&\leq \left(\int_{\mathbb{R}^{n+1}}\int_{1}^{2}|\mathcal{U}_{k,t}f(\overline{x})|^pdtd\overline{x}\right)^{1/p}\\
&=\left(\int_{1}^{2}\int_{\mathbb{R}^{n+1}}|\mathcal{U}_{k,t}f(\overline{x})|^pd\overline{x}dt\right)^{1/p}\\
&=\left(\int_{1}^{2}\int_{\mathbb{R}^{n+1}}|U_{k,t}*f(\overline{x}-tD)|^pd\overline{x}dt\right)^{1/p}\\
&\leq C\|f\|_{L^p}.
\end{align*}
Then (\ref{ineq-0910-02-00001}) follows by applying Proposition \ref{pro-001} with $B=\mathcal{O}(2^{k(n+2)})$ and $A$ being the constant just obtained.

\end{document}